\documentclass[11pt,reqno]{amsart}
\usepackage{amsmath,amssymb,mathtools}
\usepackage[a4paper,margin=1.05in]{geometry}
\usepackage[expansion=false]{microtype}
\usepackage{enumitem}
\usepackage{booktabs,array}
\usepackage{xcolor}
\usepackage{hyperref}
\hypersetup{colorlinks=true,linkcolor=blue!50!black,citecolor=green!35!black,urlcolor=blue!50!black,
  pdftitle={Categorical Torelli and Fourier--Mukai partners of hypersurfaces},
  pdfauthor={Xun Lin and Shizhuo Zhang}}
\setlist[enumerate]{label=\textup{(\roman*)},leftmargin=2.2em,itemsep=3pt}
\numberwithin{equation}{section}

\newtheorem{theorem}{Theorem}[section]
\newtheorem{proposition}[theorem]{Proposition}
\newtheorem{lemma}[theorem]{Lemma}
\newtheorem{corollary}[theorem]{Corollary}
\newtheorem{question}[theorem]{Question}
\theoremstyle{definition}
\newtheorem{definition}[theorem]{Definition}
\newtheorem{remark}[theorem]{Remark}

\newcommand{\C}{\mathbf C}
\newcommand{\Z}{\mathbf Z}
\newcommand{\Q}{\mathbf Q}
\newcommand{\PP}{\mathbf P}
\newcommand{\OO}{\mathcal O}
\newcommand{\Ku}{\mathcal{K}u}
\newcommand{\HH}{\mathrm{HH}}
\newcommand{\NHH}{\mathrm{NHH}}
\newcommand{\HS}{\mathrm{HS}}
\newcommand{\Perf}{\operatorname{Perf}}
\newcommand{\RHom}{\operatorname{RHom}}
\newcommand{\Hom}{\operatorname{Hom}}
\newcommand{\Ext}{\operatorname{Ext}}
\newcommand{\End}{\operatorname{End}}
\newcommand{\Sym}{\operatorname{Sym}}
\newcommand{\Span}{\operatorname{Span}}
\newcommand{\Spec}{\operatorname{Spec}}
\newcommand{\Proj}{\operatorname{Proj}}
\newcommand{\Cone}{\operatorname{Cone}}
\newcommand{\coker}{\operatorname{coker}}
\newcommand{\rank}{\operatorname{rank}}
\newcommand{\tr}{\operatorname{tr}}
\newcommand{\Id}{\mathrm{Id}}
\newcommand{\KS}{\mathrm{KS}}
\newcommand{\At}{\mathrm{At}}
\newcommand{\Hdg}{\operatorname{Hdg}}
\newcommand{\GL}{\operatorname{GL}}
\newcommand{\PGL}{\operatorname{PGL}}
\newcommand{\FM}{\operatorname{FM}}
\newcommand{\prim}{\mathrm{prim}}
\newcommand{\tensorL}{\otimes^{\mathrm L}}
\newcommand{\Bd}{\mathfrak B}
\renewcommand{\AA}{\mathcal A}

\title{Categorical Torelli and Fourier--Mukai partners of hypersurfaces}

\author{Xun Lin}
\address{School of Science and Engineering, The Chinese University of Hong Kong, Shenzhen, Shenzhen 518172, P.~R.~China}
\email{linxun@cuhk.edu.cn}

\author{Shizhuo Zhang}
\address{School of Mathematics, Sun Yat-sen University, Guangzhou 510275, P.~R.~China}
\email{zhangshzh28@mail.sysu.edu.cn}

\date{}

\begin{document}

\begin{abstract}
Let $X\subset\PP^{M-1}$ be a smooth hypersurface of degree $d<M$, and let $\Ku(X)$ be its Kuznetsov component, the orthogonal complement of the line bundles $\OO_X,\dots,\OO_X(M-d-1)$ in the derived category. We prove that $\Ku(X)$ determines $X$ for every smooth quartic in $\PP^{4k+1}$, $k\geq1$; together with earlier work, this gives categorical Torelli for every smooth hypersurface of degree $d<M$ with $d\nmid M$. The same algebraic reconstruction proves the classical generic Torelli theorem for these quartics, which completes generic Torelli for smooth hypersurfaces whose degree does not divide $M$, cubic surfaces excepted. When $d$ divides $M$, the Kuznetsov component is a Calabi--Yau category, as for cubic fourfolds, where categorical Torelli is known to fail. We show that, apart from cubic fourfolds, such a hypersurface has a number of Fourier--Mukai partners bounded by a constant depending only on $d$ and $M$, and that categorical Torelli holds in sufficiently large degree, in every dimension. The main tool is a deformation theorem: when geometric deformations surject onto the second Hochschild cohomology of the Kuznetsov components, an equivalence deforms together with the hypersurfaces, so statements proved for general hypersurfaces extend to all smooth ones.
\end{abstract}

\maketitle
\setcounter{tocdepth}{1}
\tableofcontents

\section{Introduction}

\subsection{Background}
Bondal and Orlov proved that a smooth projective variety with ample or anti-ample canonical bundle is determined by its bounded derived category \cite{BO01}. When the anticanonical bundle is ample, one can often say more. The derived category then usually contains an exceptional collection, and its orthogonal complement, the \emph{Kuznetsov component}, is a much smaller category. One says that a family of such varieties satisfies \emph{categorical Torelli} if the Kuznetsov component determines the variety. The name comes from the classical Torelli theorems, which recover a variety from its Hodge structure.

Smooth hypersurfaces are the basic examples. Let $X\subset\PP^{M-1}_\C$ be a smooth hypersurface of degree $3\leq d<M$. The line bundles $\OO_X,\OO_X(1),\dots,\OO_X(M-d-1)$ form an exceptional collection, and the Kuznetsov component is
\begin{equation}\label{eq:ku}
 \Ku(X)=\langle\OO_X,\OO_X(1),\dots,\OO_X(M-d-1)\rangle^\perp.
\end{equation}
It carries a natural dg enhancement. Throughout the paper, equivalences are $\C$-linear equivalences of these enhancements, given by perfect Fourier--Mukai kernels; they are not required to commute with the rotation functor of $\Ku(X)$, the autoequivalence induced by $\OO_X(1)$.

The answer to the reconstruction question depends on whether $d$ divides $M$. The Kuznetsov component has a rotation functor $\rho$ with $\rho^d\simeq[2]$, and its Serre functor is $\rho^{-M}[M]$ \cite{Kuz19,LRZ26}. If $d\nmid M$, the Serre functor retains part of the information carried by $\rho$. If $M=md$, the Serre functor is the shift $[m(d-2)]$, so $\Ku(X)$ is a \emph{Calabi--Yau category} of dimension $m(d-2)$. For a cubic fourfold it is a K3 category \cite{Kuz10}. In the Calabi--Yau case the Serre functor says nothing about $\rho$, and nonisomorphic hypersurfaces can have equivalent Kuznetsov components. It is then natural to count them. For fixed $d$ and $M$ we put
\begin{equation}\label{eq:fm}
 \FM_{d,M}(X)=\{Y\subset\PP^{M-1}\ \text{smooth of degree }d:\ \Ku(Y)\simeq\Ku(X)\}/\simeq_{\mathrm{proj}},
\end{equation}
where $\simeq_{\mathrm{proj}}$ is projective isomorphism, and we call the elements of $\FM_{d,M}(X)$ the \emph{Fourier--Mukai partners} of $X$.\footnote{For hypersurfaces, a Fourier--Mukai partner of $X$ thus means a smooth hypersurface $Y$ of the same degree and dimension whose Kuznetsov component is equivalent to that of $X$; the full derived categories are not required to be equivalent. Indeed $\mathrm D^{\mathrm b}(X)\simeq\mathrm D^{\mathrm b}(Y)$ already forces $X\simeq Y$ by \cite{BO01}, since $d<M$. For cubic fourfolds this is the notion used in \cite{Huy17,Per21}.} This set contains $X$, and categorical Torelli for $X$ means that it has one element.

\subsection{Previous results}
Categorical Torelli theorems are known for many varieties with ample anticanonical bundle, starting with cubic threefolds \cite{BMMS12}. For hypersurfaces, Huybrechts and Rennemo showed that the Hochschild cohomology of $\Ku(X)$, together with the action of the rotation functor, recovers the Jacobian ring; this gives categorical Torelli for the pair $(\Ku(X),\rho)$, in particular for cubic fourfolds \cite{HR19}. Without the rotation functor, Pirozhkov \cite{Pir24} and Lin and Zhang \cite{LZ25} proved categorical Torelli for hypersurfaces of degree $d<M$ under arithmetic conditions on $d$ and $M$. Lin, Rennemo and Zhang proved it whenever $d\nmid M$, except for quartics in $\PP^{4k+1}$ \cite[Theorem~3.12]{LRZ26}. For these quartics they proved a generic version when $k\geq150$ \cite[Corollary~3.13]{LRZ26}, using Voisin's generic Torelli theorem \cite{Voi22}. The present paper completes the case $d\nmid M$ and treats the case $d\mid M$; the resulting picture is summarized in Subsection~\ref{subsec:summary}.

The quartics in $\PP^{4k+1}$ are also the last open case of a classical problem. Generic Torelli for hypersurfaces asks whether a very general hypersurface is determined by the Hodge structure on its primitive middle cohomology. Griffiths described this Hodge structure through the Jacobian ring \cite{Gri69}, and Carlson and Griffiths showed that its infinitesimal variation is given by multiplication in that ring \cite{CG80}. Building on this, Donagi proved generic Torelli for smooth hypersurfaces of degree $d\geq3$ in $\PP^{M-1}$, $M\geq4$, except for cubic surfaces, the case $d\mid M$, quartics in $\PP^{4k+1}$, and sextics in $\PP^{6k+2}$ \cite{Don83}; see also \cite{DG84,CDT87}. Cox and Green removed the sextic exception \cite{CG90}, and Voisin treated the case $d\mid M$ in sufficiently large degree \cite{Voi22}. For the quartic series, Voisin's methods recover the polynomial structure in large dimension, which is the input behind \cite[Corollary~3.13]{LRZ26}.

\subsection{Main results}\label{subsec:main}

Our first result settles the remaining series of the case $d\nmid M$.

\begin{theorem}\label{thm:quartic}
For every $k\geq1$, let $X,Y\subset\PP^{4k+1}_\C$ be smooth quartic hypersurfaces. An enhanced Fourier--Mukai equivalence
\[
 \Ku(X)\simeq\Ku(Y)
\]
implies that $X$ and $Y$ are projectively isomorphic.
\end{theorem}

Together with \cite[Theorem~3.12]{LRZ26}, this gives categorical Torelli for every smooth hypersurface of degree $d<M$ with $d\nmid M$ (Corollary~\ref{cor:nondivisible}). The proof goes through a generic reconstruction theorem for all $k\geq1$ (Theorem~\ref{thm:generic-quartic}): a general quartic in $\PP^{4k+1}$ is determined by the even part of its Jacobian ring.

The same reconstruction, combined with Donagi's symmetrizer construction and the theory of Hodge loci, settles the classical generic Torelli problem for these quartics. Together with the work of Donagi and of Cox and Green, this gives the following theorem, proved in Appendix~\ref{app:generic}; here the degree may also be at least $M$.

\begin{theorem}\label{thm:classical}
Let $d\geq3$ and $M\geq4$, with $d\nmid M$ and $(d,M)\neq(3,4)$. If $X\subset\PP^{M-1}$ is a very general smooth hypersurface of degree $d$ and $Y\subset\PP^{M-1}$ is any smooth hypersurface of degree $d$ with
\[
 H^{M-2}_{\prim}(X,\Q)\simeq H^{M-2}_{\prim}(Y,\Q)
\]
as rational Hodge structures, then $X$ and $Y$ are projectively isomorphic.
\end{theorem}

Cubic surfaces are excluded because their primitive cohomology is entirely of type $(1,1)$. The new case in Theorem~\ref{thm:classical} is that of quartics in $\PP^{4k+1}$ (Theorem~\ref{thm:quartic-generic-torelli}).

When $d$ divides $M$, the Kuznetsov component is a Calabi--Yau category. The best-known example is the cubic fourfold, $(d,M)=(3,6)$, whose Kuznetsov component is a K3 category \cite{Kuz10}, and here the classical and the categorical points of view give strikingly different answers. Classically, a cubic fourfold is determined by its Hodge structure: this is Voisin's Torelli theorem \cite{Voi86}. Hassett showed that the Hodge structures of certain special cubic fourfolds contain that of a K3 surface \cite{Has00}, and Kuznetsov showed that for some of them $\Ku(X)$ is the derived category of a K3 surface \cite{Kuz10}. Since a K3 surface can have nonisomorphic Fourier--Mukai partners, this suggests that $\Ku(X)$ alone need not determine $X$, and this is what happens. Huybrechts and Rennemo proved that $\Ku(X)$ together with its rotation functor determines $X$ \cite{HR19}. Without the rotation functor, Huybrechts proved that $\Ku(X)$ determines a very general cubic fourfold, and that every cubic fourfold has only finitely many Fourier--Mukai partners \cite{Huy17}, just as for K3 surfaces and, more generally, for surfaces \cite{BM01}. Pertusi then found very general special cubic fourfolds with nonisomorphic Fourier--Mukai partners and showed that their number can be arbitrarily large \cite{Per21}, and Fan and Lai computed this number for every very general special cubic fourfold \cite{FL25}. So for cubic fourfolds categorical Torelli fails, but in a controlled way: partners exist, and each cubic fourfold has only finitely many of them.

Every hypersurface whose degree divides $M$ has a Calabi--Yau Kuznetsov component, so the question raised by cubic fourfolds makes sense for all of them: \emph{does a smooth hypersurface have only finitely many Fourier--Mukai partners among hypersurfaces of the same degree and dimension?} To state our answer we fix
\begin{equation}\label{eq:cy-range}
 d\geq3,\qquad m\geq2,\qquad M=md,\qquad (d,m)\neq(3,2),
\end{equation}
so that cubic fourfolds are excluded, and write $c=m(d-2)\geq3$ for the Calabi--Yau dimension of $\Ku(X)$. The answer is positive, and in a stronger form than is possible for cubic fourfolds: the number of partners is bounded independently of the hypersurface.

\begin{theorem}\label{thm:finite}
There is an integer $\Bd_{d,m}$, depending only on $d$ and $m$, such that every smooth hypersurface $X\subset\PP^{md-1}$ of degree $d$ satisfying~\eqref{eq:cy-range} has
\[
 1\leq\#\FM_{d,md}(X)\leq\Bd_{d,m}.
\]
\end{theorem}

The bound comes from an intrinsic graded algebra attached to $\Ku(X)$, the subalgebra of Hochschild cohomology generated by $\HH^2$: we bound the number of hypersurfaces whose algebras of this kind are isomorphic. Cubic fourfolds are excluded because there the second Hochschild cohomology of $\Ku(X)$ is two dimensions larger than the space of deformations of $X$ (Remark~\ref{rem:cubic-fourfold}).

Theorem~\ref{thm:finite} does not say whether partners actually exist. Since categorical Torelli fails for cubic fourfolds, one might expect it to fail for every hypersurface with a Calabi--Yau Kuznetsov component: there the Serre functor is a shift, so the Hochschild--Serre algebra carries no information about the rotation functor. Our next result shows that the opposite happens in large degree.

\begin{theorem}\label{thm:large}
There is an integer $d_0\geq5$, independent of the ambient dimension, such that for every $d\geq d_0$ and every $M>d$, two smooth degree-$d$ hypersurfaces $X,Y\subset\PP^{M-1}$ with equivalent enhanced Kuznetsov components are projectively isomorphic. Equivalently, $\#\FM_{d,M}(X)=1$.
\end{theorem}

For $d\mid M$ this follows from Voisin's rational generic Torelli theorem \cite[Theorem~0.7(1)]{Voi22} and Theorem~\ref{thm:deformation}; the threshold is uniform in $M$ by \cite[Theorem~0.11]{Voi22}. For $d\nmid M$ it follows from \cite[Theorem~3.12]{LRZ26}. So the failure of categorical Torelli for cubic fourfolds is a phenomenon of low degree.

\begin{remark}\label{rem:wu}
Theorem~\ref{thm:large} says nothing about low degree, but categorical Torelli can still hold for special hypersurfaces there. As pointed out to us by Dongjian Wu, it holds for a subfamily of the cubic hypersurfaces of dimension $3m-2$ with $m\geq3$, whose Kuznetsov components are Calabi--Yau of dimension $m$: if a smooth cubic $X\subset\PP^{3m-1}$ is, after a linear change of coordinates, defined by an equation of the form
\[
 f(x_0,\dots,x_{a-1})+g(x_a,\dots,x_{3m-1})=0,\qquad 1\leq a\leq3m-1,
\]
that is, a Sebastiani--Thom sum of two cubic forms in disjoint sets of variables, then every smooth cubic $Y\subset\PP^{3m-1}$ with $\Ku(X)\simeq\Ku(Y)$ is projectively isomorphic to $X$. Here $Y$ is not assumed to be of this form, and the equivalence need not commute with the rotation functor.
\end{remark}

To see where the boundary between low and large degree comes from, note that the proof of Theorem~\ref{thm:large} for $d\mid M$ has two steps. The first step shows that $\Ku(X)$ determines $X$ when $X$ is very general: an equivalence preserves topological K-theory with its Hodge structure, and for very general $X$ this recovers the rational Hodge structure on $H^{M-2}_{\prim}(X,\Q)$ (Proposition~\ref{prop:hodge}), so a generic Torelli theorem for this Hodge structure gives uniqueness. The second step extends uniqueness from very general hypersurfaces to all of them by deforming the equivalence (Theorem~\ref{thm:deformation} below), which requires every first-order deformation of $\Ku(X)$ to come from a deformation of $X$.

For cubic fourfolds the first step works \cite{Huy17}, but the second fails. The second Hochschild cohomology of $\Ku(X)$ is two dimensions larger than the space of deformations of $X$, as for a K3 surface, so an equivalence between two special cubic fourfolds need not extend to nearby cubic fourfolds; and indeed Fourier--Mukai partners appear on Hassett's divisors of special cubic fourfolds \cite{Huy17,Per21}. For every other hypersurface with $d\mid M$ the second step works, because the additional Hochschild classes of $\Ku(X)$ lie in degree $c=m(d-2)\geq3$ (Proposition~\ref{prop:veronese-intrinsic}), and in low degree it is the first step that is missing. The reason is that, without the rotation functor, the algebra we can extract from $\Ku(X)$ is the $d$-th Veronese $R_f^{(d)}=\bigoplus_a(R_f)_{ad}$ of the Jacobian ring, and the classical side faces the same limitation. By Griffiths' residue description the Hodge pieces of $X$ are graded pieces $(R_f)_{(q+1)d-M}$ of the Jacobian ring, all of degree divisible by $d$ when $d\mid M$, so the infinitesimal variation of Hodge structure also sees only $R_f^{(d)}$. This is why the case $d\mid M$ is excluded from Donagi's theorem; Voisin treated it in large degree, using higher-order variations of Hodge structure along Schiffer variations \cite{Voi22}. In low degree the following question remains.

\begin{question}\label{q:generic-torelli}
Let $d\geq3$ and $m\geq2$ with $(d,m)\neq(3,2)$, and put $M=md$. Is a very general smooth hypersurface $X\subset\PP^{M-1}$ of degree $d$ determined, up to projective isomorphism, among all smooth hypersurfaces of degree $d$ in $\PP^{M-1}$, by the rational Hodge structure $H^{M-2}_{\prim}(X,\Q)$?
\end{question}

By \cite{Voi22} the answer is positive for $d$ sufficiently large; as far as we know, the question is open in low degree, for instance for quartic sixfolds in $\PP^7$ and cubic sevenfolds in $\PP^8$. A positive answer would also settle categorical Torelli for these hypersurfaces. Since the second step works for them, Proposition~\ref{prop:conditional} shows that if Question~\ref{q:generic-torelli} has a positive answer for some $(d,m)$, then any two smooth hypersurfaces of degree $d$ in $\PP^{md-1}$ with equivalent Kuznetsov components are projectively isomorphic. Thus, apart from cubic fourfolds, categorical Torelli for hypersurfaces with a Calabi--Yau Kuznetsov component reduces to a classical question about Hodge structures, and Theorem~\ref{thm:finite} gives its finite version unconditionally.

Theorems~\ref{thm:quartic}, \ref{thm:finite} and~\ref{thm:large} rest on a deformation argument. In each case the reconstruction works directly only for a general hypersurface: for quartics because the genericity condition used in Section~\ref{sec:quartics} fails, for instance, for the Fermat quartic; in the Calabi--Yau case because the count is made on a slice near the Fermat hypersurface; and in large degree because Hodge theory is used at a very general point. To pass from a general hypersurface to all of them, we deform the equivalence together with the hypersurfaces. Suppose that $\Ku(X)\simeq\Ku(Y)$, and let $X$ move in a family. First-order deformations of a category are measured by its second Hochschild cohomology, so the equivalence can follow a deformation of $X$ only if the induced deformation of $\Ku(X)$ is matched by a deformation of $\Ku(Y)$. If every class in $\HH^2(\Ku(Y))$ comes from a deformation of $Y$, such a matching deformation of $Y$ always exists. The equivalence then extends over a whole family of pairs, and statements proved for general members spread to all members.

The following theorem makes this precise for arbitrary families. For an admissible component $\AA_X\subset\Perf(X)$, let
\begin{equation}\label{eq:eta}
 \eta_X:H^1(X,T_X)\longrightarrow\HH^2(X)\longrightarrow\HH^2(\AA_X)
\end{equation}
be the geometric summand of the Hochschild--Kostant--Rosenberg map followed by restriction; it sends a first-order deformation of $X$ to the induced first-order deformation of $\AA_X$. For a smooth family $\mathcal X\to B$ with a relative component, let
\begin{equation}\label{eq:kappa}
 \kappa_b:T_bB\xrightarrow{\ \KS_b\ }H^1(\mathcal X_b,T_{\mathcal X_b})\xrightarrow{\ \eta_{\mathcal X_b}\ }\HH^2(\AA_b)
\end{equation}
be its composite with the Kodaira--Spencer map, so that $\kappa_b$ records which deformations of the component are realized inside the family. By an \emph{individual very-general locus} in an irreducible complex variety we mean the complement of a countable union of proper closed algebraic subsets.

\begin{theorem}\label{thm:deformation}
Let $\mathcal X_i\to B_i$ be smooth projective polarized families over smooth irreducible complex varieties, carrying relative admissible components whose left and right projection functors have perfect kernels compatible with base change.
\begin{enumerate}
\item Fix $b_i\in B_i(\C)$ and an equivalence $F:\AA_{1,b_1}\simeq\AA_{2,b_2}$. If the maps $\kappa_{b_i}$ are surjective and $\HH^j(\AA_{i,b_i})=0$ for $j<0$, then there are a smooth irreducible variety $T$, a point $t\in T$, smooth morphisms $p_i:T\to B_i$ with $p_i(t)=b_i$, and a family of equivalence kernels over $T$ restricting to $F$ at $t$.
\item Let $\mathcal X\to B$ be one such family, and suppose that these conditions hold at every fibre. Assume that a morphism $B\to\mathcal M$ to a separated algebraic space distinguishes polarized isomorphism classes. If an integer $N$ bounds the size of every set of pairwise nonisomorphic fibres over an individual very-general locus whose components are mutually equivalent, then the same bound holds over all of $B$.
\end{enumerate}
In particular, under these hypotheses, generic categorical Torelli implies categorical Torelli for every smooth member of the family.
\end{theorem}

Part~(i) says that, when the maps $\kappa$ are surjective and negative Hochschild cohomology vanishes, an equivalence between the components of two fibres extends to a family of equivalences over a variety $T$ that maps smoothly, and hence dominantly, onto both parameter spaces. Part~(ii) is the consequence used in this paper: a bound on the number of partners that holds over a general locus holds everywhere. To apply Theorem~\ref{thm:deformation} to hypersurfaces, we prove in Theorem~\ref{thm:deg4} and Proposition~\ref{prop:cubic-def} that, for smooth hypersurfaces of degree $4\leq d<M$ and for cubics in $\PP^{3m-1}$ with $m\geq3$,
\begin{equation}\label{eq:hkr-iso}
 H^1(X,T_X)\xrightarrow{\ \sim\ }\HH^2(\Ku(X)),\qquad\HH^{<0}(\Ku(X))=0.
\end{equation}

\subsection{Categorical Torelli for hypersurfaces: the current picture}\label{subsec:summary}
Combining the results of this paper with earlier work gives the following picture for smooth hypersurfaces $X\subset\PP^{M-1}$ of degree $3\leq d<M$.
\begin{center}
\renewcommand{\arraystretch}{1.2}
\begin{tabular}{>{\raggedright\arraybackslash}p{0.38\textwidth}>{\raggedright\arraybackslash}p{0.3\textwidth}>{\raggedright\arraybackslash}p{0.22\textwidth}}
\toprule
Case & Categorical Torelli & Reference\\
\midrule
$d\nmid M$, except quartics in $\PP^{4k+1}$ & holds & \cite{Pir24,LZ25,LRZ26}\\
$d=4$, $M=4k+2$, $k\geq1$ & holds & Theorem~\ref{thm:quartic}\\
$d\geq d_0$, any $M>d$ & holds & Theorem~\ref{thm:large}\\
cubic fourfolds ($d=3$, $M=6$) & fails in general, with finitely many partners; holds for very general $X$, and with the rotation functor & \cite{Huy17,Per21}; \cite{HR19}\\
$d\mid M$, $d<d_0$, not a cubic fourfold & finitely many partners, uniformly bounded; holds if Question~\ref{q:generic-torelli} has a positive answer & Theorem~\ref{thm:finite}; Proposition~\ref{prop:conditional}\\
\bottomrule
\end{tabular}
\end{center}
On the classical side, Theorem~\ref{thm:classical} completes generic Torelli for smooth hypersurfaces whose degree does not divide $M$, cubic surfaces excepted. On the categorical side, categorical Torelli now holds for \emph{every} smooth hypersurface of degree $d<M$ with $d\nmid M$, which completes the line of work in \cite{Pir24,LZ25,LRZ26}, and for every smooth hypersurface with $d_0\leq d<M$, in every dimension. When $d$ divides $M$, the Kuznetsov component is a Calabi--Yau category and the situation is different: for cubic fourfolds the Kuznetsov component alone does not always determine $X$. The question remains open only for hypersurfaces with $d\mid M$ and $d<d_0$ other than cubic fourfolds, for instance quartics in $\PP^7$ and cubic sevenfolds; for these, Theorem~\ref{thm:finite} shows that each hypersurface has finitely many partners, with a bound depending only on $d$ and $M$, and Proposition~\ref{prop:conditional} reduces categorical Torelli to Question~\ref{q:generic-torelli}.

\subsection{Outline of the proofs}
We first explain the proof of Theorem~\ref{thm:deformation}, on which the other results depend. Let $P$ be the kernel of an equivalence between two admissible components, extended to the ambient varieties by the projection functors. By derived Morita theory, $\Ext^2(P,P)$ is the second Hochschild cohomology of either component. Over a small extension of Artinian rings, we fix a lift of the first variety and choose a lift of the second. The obstruction to deforming $P$ lies in $\Ext^2(P,P)$, and changing the lift of the second variety changes the obstruction by the image of the change under $\kappa$. If $\kappa$ is surjective, the lift of the second variety can be chosen so that the obstruction vanishes. Projection to the relative components and derived Nakayama show that the lifted kernel is again an equivalence. Using the algebraic stack of complexes, we obtain a variety $T$ with smooth maps to both bases carrying a family of equivalences. Taking fibre products of such varieties deforms several equivalences at once, which gives part~(ii).

To apply Theorem~\ref{thm:deformation} to the universal family of hypersurfaces, we must show that every first-order deformation of $\Ku(X)$ comes from a deformation of $X$. For $d\geq4$, the exceptional collection in~\eqref{eq:ku} has pseudoheight $d-1$, so by Kuznetsov's theory of normal Hochschild cohomology the restriction $\HH^2(X)\to\HH^2(\Ku(X))$ is injective. The normal sequence identifies $H^1(X,T_X)$ with the degree-$d$ part of the Jacobian ring, and the matrix-factorization description of $\Ku(X)$ gives a target of the same dimension. For Calabi--Yau cubics of dimension at least seven, injectivity comes instead from the action on the two highest primitive Hodge pieces.

With this in place, Theorem~\ref{thm:quartic} reduces to a statement about general quartics. For quartics in $\PP^{4k+1}$, the functor $S_{\Ku(X)}^{-1}[2k+2]$ is the square of the rotation functor, and the natural transformations from the identity to its powers recover the even part $\bigoplus_a(R_f)_{2a}$ of the Jacobian ring $R_f$. Its degree-one part is the space $S_2$ of quadrics. The quadratic relations of the Veronese embedding are exactly the relations of rank at most three among products of quadrics, and we show that a general quartic introduces no further relations of rank at most three. This lets us recover the polynomial structure, then the Jacobian ideal, and then the quartic by Donagi's reconstruction theorem. A uniform dimension estimate treats $\PP^{4k+1}$ with $k\geq3$, and separate calculations treat $k=1,2$. Theorem~\ref{thm:deformation} then extends the generic statement to all smooth quartics.

The classical Theorem~\ref{thm:classical} uses the same algebraic reconstruction with a different input. The Hodge structure at a single point does not see the multiplication of the Jacobian ring, but its infinitesimal variation does, and Donagi's symmetrizer construction recovers the even Jacobian algebra from it. Algebraic Hodge loci then show that, for a very general quartic, an isomorphism of Hodge structures comes from an isomorphism of infinitesimal variations. Combined with the results of Donagi and of Cox and Green, this gives the theorem.

Theorem~\ref{thm:finite} follows the same pattern, with a count in place of a reconstruction. When $d\mid M$, the subalgebra of Hochschild cohomology generated by $\HH^2$ is the $d$-th Veronese $\bigoplus_a(R_f)_{ad}$ of the Jacobian ring, an intrinsic invariant of $\Ku(X)$, and we bound the number of hypersurfaces with isomorphic Veronese algebras. Near the Fermat hypersurface we choose a slice transverse to the action of $\GL_M$, compute the infinitesimal stabilizer of the top form of the algebra, and show that the map from the slice and the changes of basis to the space of top forms has finite fibres. This map is quasi-finite, so its fibres have bounded size, which gives a bound depending only on $d$ and $m$. Theorem~\ref{thm:deformation} carries the bound from the slice to every smooth hypersurface.

Finally, for Theorem~\ref{thm:large} the generic input comes from Hodge theory. An equivalence preserves topological K-theory with its Hodge structure and Euler form. For a very general hypersurface this recovers the primitive rational Hodge structure, and Voisin's generic Torelli theorem gives uniqueness. Theorem~\ref{thm:deformation} extends uniqueness to every smooth hypersurface.

\subsection{Organization}
Section~\ref{sec:deformation-thm} proves Theorem~\ref{thm:deformation}. Section~\ref{sec:hochschild} computes the deformation maps and the relevant Hochschild cohomology of Kuznetsov components of hypersurfaces. Section~\ref{sec:quartics} proves Theorem~\ref{thm:quartic}, Section~\ref{sec:finiteness} proves Theorem~\ref{thm:finite}, and Section~\ref{sec:large} proves Theorem~\ref{thm:large}, together with Proposition~\ref{prop:conditional}, which reduces the low-degree Calabi--Yau cases to Question~\ref{q:generic-torelli}. Appendix~\ref{app:generic} applies the algebraic reconstruction of Section~\ref{sec:quartics} to the classical problem and proves Theorem~\ref{thm:classical}.

\subsection*{Acknowledgements}
Shizhuo Zhang is supported by the Sun Yat-sen University Starting Grant 34000-12256019. Xun Lin is supported by the Chinese University of Hong Kong(Shen Zhen) starting grant.

\subsection*{AI disclosure}
The results of this paper were obtained through interaction with ChatGPT~6, and the authors verified all results and proofs independently. The paper was written by the authors, with the assistance of AI tools in editing the exposition. The authors take full responsibility for its content.

\section{Deformation of equivalences}\label{sec:deformation-thm}

In this section we prove Theorem~\ref{thm:deformation}. Subsection~\ref{sec:kernels} identifies the self-extensions of an equivalence kernel with Hochschild cohomology, Subsection~\ref{sec:deform} lifts equivalence kernels through small extensions, and Subsection~\ref{sec:special} turns this infinitesimal statement into algebraic correspondences and proves the specialization results.

\subsection{Projected kernels and Hochschild cohomology}\label{sec:kernels}

This subsection fixes the relative setting and identifies the self-extensions of an equivalence kernel with Hochschild cohomology. All schemes are of finite type over $\C$, except for the spectra of local Artinian $\C$-algebras used in deformation arguments. All categories carry their given dg enhancements. Tensor products of categories are derived, and their perfect hulls are understood.

For $i=1,2$, let
\[
 \pi_i:\mathcal X_i\longrightarrow B_i
\]
be a smooth projective morphism with $B_i$ smooth. We fix a $B_i$-linear admissible component $\AA_i\subset\Perf(\mathcal X_i)$, and assume that its left and right projection functors are represented by perfect relative kernels and commute with base change. The same convention applies after base change to a local Artinian $\C$-algebra. We write $\AA_{i,S}$ for the component after a morphism $S\to B_i$, and $j_{i,S}$, $r_{i,S}$ for its inclusion and the right adjoint of the inclusion. These hypotheses also give relative projection kernels on products, by applying the two factorwise projections in turn. Their compatibility with base change is used both for infinitesimal lifting and for the construction of an open stack of equivalence kernels.

Fix $b_i\in B_i(\C)$ and set $X_i=\mathcal X_{i,b_i}$ and $\AA_{i,0}=\AA_{i,b_i}$; write $j_i$ and $r_i$ for the inclusion and its right adjoint at this fibre. An equivalence
\[
 F:\AA_{1,0}\xrightarrow{\ \sim\ }\AA_{2,0}
\]
always means an enhanced $\C$-linear equivalence. Its \emph{projected kernel} is the perfect Fourier--Mukai kernel $P\in\Perf(X_1\times X_2)$ representing
\begin{equation}\label{eq:projected}
 j_2Fr_1:\Perf(X_1)\longrightarrow\Perf(X_2).
\end{equation}
For a smooth proper dg category $\AA$ we use the convention
\[
 \HH^m(\AA)=\Hom_{\Perf(\AA^{\mathrm{op}}\otimes\AA)}(\Id_\AA,\Id_\AA[m]),
\]
where the identity functor is the diagonal bimodule. By derived Morita theory, its derived endomorphisms compute the Hochschild cochain complex \cite[Corollary~8.1]{Toe07}. For smooth projective varieties, enhanced functors are described by geometric kernels \cite[Theorem~8.15]{Toe07}.

\begin{lemma}\label{lem:ext-hh}
Let $P$ be the projected kernel of $F$. For every integer $m$ there are isomorphisms
\begin{equation}\label{eq:ext-hh}
 \HH^m(\AA_{1,0})\xrightarrow{\ \chi_{P,1}\ }\Ext^m_{X_1\times X_2}(P,P)\xleftarrow{\ \chi_{P,2}\ }\HH^m(\AA_{2,0}),
\end{equation}
given by the characteristic actions on the source and the target of the equivalence bimodule. The characteristic action of $\HH^m(X_i)$ on $P$ factors through the restriction map $\rho_i:\HH^m(X_i)\to\HH^m(\AA_{i,0})$.
\end{lemma}

\begin{proof}
Let $\AA_{1,0}^\vee\subset\Perf(X_1)$ be the image of $\AA_{1,0}$ under derived duality. The exterior product realizes the category of perfect $\AA_{1,0}$--$\AA_{2,0}$ bimodules as
\begin{equation}\label{eq:bimod-realization}
 \AA_{1,0}^\vee\boxtimes\AA_{2,0}\subset\Perf(X_1\times X_2).
\end{equation}
Indeed, for $E,E'\in\AA_{1,0}$ and $G,G'\in\AA_{2,0}$, duality and the K\"unneth formula give
\[
 \RHom_{X_1\times X_2}\bigl(E^\vee\boxtimes G,(E')^\vee\boxtimes G'\bigr)\simeq\RHom_{X_1}(E',E)\otimes_\C\RHom_{X_2}(G,G'),
\]
which proves full faithfulness on the exterior products, and hence on their perfect hull. The kernel $E^\vee\boxtimes G$ sends $V\in\Perf(X_1)$ to $\RHom(E,V)\otimes G\simeq\RHom(E,r_1V)\otimes G$, so the realization~\eqref{eq:bimod-realization} gives exactly the extension~\eqref{eq:projected}.

The bimodule representing $F$ is invertible. Tensoring with its inverse identifies its derived endomorphism complex with that of the diagonal bimodule of $\AA_{1,0}$, and tensoring on the other side gives the corresponding identification for $\AA_{2,0}$. By \cite[Corollary~8.1]{Toe07}, these complexes compute the Hochschild cochain complexes, and full faithfulness of~\eqref{eq:bimod-realization} gives~\eqref{eq:ext-hh}. The maps are the characteristic actions, since they are induced by composing an endomorphism of an identity bimodule with the equivalence bimodule; the same description for a projection kernel is given in \cite[Theorem~4.5]{Kuz09}.

For the last assertion, represent a class in $\HH^m(X_i)$ by an endomorphism of the diagonal kernel. Its action on $P$ is convolution on the corresponding factor. Since $P$ lies in~\eqref{eq:bimod-realization}, inserting the projection to $\AA_{i,0}$ before this convolution does not change the action, and the induced endomorphism of the identity bimodule of $\AA_{i,0}$ is $\rho_i$ of the original class; this is the restriction construction of \cite[Sections~4 and~7]{Kuz09}. Hence the ambient action is $\chi_{P,i}\circ\rho_i$.
\end{proof}

For a smooth projective variety $X$ with an admissible component $\AA_X$, recall the map $\eta_X$ of~\eqref{eq:eta}. Its first arrow is the geometric summand of the Hochschild--Kostant--Rosenberg decomposition; equivalently, a Kodaira--Spencer class acts by contraction with the Atiyah class of the diagonal. For the families $\pi_i$ put
\begin{equation}\label{eq:kappa-i}
 \kappa_i:T_{b_i}B_i\xrightarrow{\ \KS_i\ }H^1(X_i,T_{X_i})\xrightarrow{\ \eta_{X_i}\ }\HH^2(\AA_{i,0}).
\end{equation}
If the Kodaira--Spencer map is surjective, surjectivity of $\eta_{X_i}$ implies surjectivity of $\kappa_i$.

\subsection{Deformation of equivalence kernels}\label{sec:deform}

We keep the relative hypotheses of Subsection~\ref{sec:kernels}. The main result of this subsection, Theorem~\ref{thm:lift}, lifts an equivalence kernel through every small extension, provided the maps $\kappa_i$ are surjective. A \emph{small extension} is an exact sequence of local Artinian $\C$-algebras with residue field $\C$,
\begin{equation}\label{eq:small}
 0\longrightarrow J\longrightarrow R'\longrightarrow R\longrightarrow0,\qquad\mathfrak m_{R'}J=0.
\end{equation}
In particular $J$ is a finite-dimensional vector space and $J^2=0$. A \emph{deformation} of $(b_1,b_2,P)$ over $R$ consists of morphisms $b_{i,R}:\Spec R\to B_i$ extending $b_i$, a projected perfect kernel $P_R$ inducing an equivalence $F_R:\AA_{1,R}\simeq\AA_{2,R}$, and an identification of its closed fibre with $P$.

\begin{theorem}\label{thm:lift}
Assume that both maps $\kappa_i$ in~\eqref{eq:kappa-i} are surjective. Let~\eqref{eq:small} be a small extension and $(b_{1,R},b_{2,R},P_R)$ a deformation of $(b_1,b_2,P)$ over $R$. Every lift $b_{1,R'}$ of $b_{1,R}$ can be completed to a deformation $(b_{1,R'},b_{2,R'},P_{R'})$ over $R'$, where $P_{R'}$ is again a projected equivalence kernel. The same holds with the two factors interchanged. Consequently, an equivalence kernel lifts through every nilpotent extension of local Artinian $\C$-algebras, with either geometric deformation prescribed.
\end{theorem}

\begin{proof}
\emph{Step 1: the obstruction group.} Since $B_2$ is smooth, choose a lift $b^0_{2,R'}$ of $b_{2,R}$, and write
\[
 Z_R=\mathcal X_{1,R}\times_R\mathcal X_{2,R},\qquad Z^0_{R'}=\mathcal X_{1,R'}\times_{R'}\mathcal X^0_{2,R'}.
\]
The square-zero immersion $Z_R\subset Z^0_{R'}$ has ideal $\OO_{Z_R}\otimes_RJ$. The obstruction to lifting $P_R$ to a perfect complex on $Z^0_{R'}$ lies in
\begin{equation}\label{eq:obs-group}
 \Ext^2_{Z_R}\bigl(P_R,P_R\tensorL_RJ\bigr),
\end{equation}
and a lift exists if and only if it vanishes. This is the obstruction theory for perfect complexes of \cite[Theorem~3.1.1]{Lie06}, with the Atiyah-class description of \cite[Corollary~3.4 and Section~3.4]{HT10}, in the corrected form of \cite{HT14}. Let $\iota:X_1\times X_2\hookrightarrow Z_R$ be the closed fibre. Since $J$ is annihilated by the maximal ideal, $P_R\tensorL_RJ\simeq\iota_*(P\otimes_\C J)$, and adjunction identifies~\eqref{eq:obs-group} with
\begin{equation}\label{eq:obs-closed}
 \Ext^2_{X_1\times X_2}(P,P)\otimes_\C J;
\end{equation}
compare \cite[Corollary~3.1.2]{Lie06}.

\smallskip\noindent\emph{Step 2: dependence on the second lift.} Let $o(P_R;b_{1,R'},b^0_{2,R'})$ be the obstruction in~\eqref{eq:obs-closed}. The lifts of $b_{2,R}$ form a torsor under $T_{b_2}B_2\otimes J$, and for $\xi$ in this space
\begin{equation}\label{eq:obs-change}
 o(P_R;b_{1,R'},b^0_{2,R'}+\xi)-o(P_R;b_{1,R'},b^0_{2,R'})=(\chi_{P,2}\otimes\Id_J)\bigl((\kappa_2\otimes\Id_J)(\xi)\bigr),
\end{equation}
with the torsor convention chosen so that the sign is as displayed. To prove~\eqref{eq:obs-change}, write the obstruction as the composition
\begin{equation}\label{eq:obs-atiyah}
 P_R\xrightarrow{\ \At(P_R)\ }P_R\tensorL\mathbf L_{Z_R/\C}[1]\xrightarrow{\ \Id\otimes\KS\ }P_R\tensorL_RJ[2],
\end{equation}
where $\mathbf L_{Z_R/\C}$ is the cotangent complex; its truncation may be used. The schemes involved are projective over Artinian $\C$-algebras, hence embed in smooth $\C$-schemes, as required in \cite{HT10,HT14}. The difference between $b^0_{2,R'}$ and $b^0_{2,R'}+\xi$ is a derivation to $J$, and functoriality of the cotangent complex identifies the resulting difference of extension classes in~\eqref{eq:obs-atiyah} with the Kodaira--Spencer class $\KS_2(\xi)$ on the second factor. Since $\mathfrak m_{R'}J=0$, its action is computed on the closed fibre $P$. Contracting this class with the Atiyah class of the second factor is its geometric Hochschild characteristic action, which Lemma~\ref{lem:ext-hh} identifies with $\chi_{P,2}\kappa_2(\xi)$. This proves~\eqref{eq:obs-change}, for an arbitrary small extension.

\smallskip\noindent\emph{Step 3: cancelling the obstruction.} By Lemma~\ref{lem:ext-hh}, $\chi_{P,2}$ is an isomorphism, and since $\kappa_2$ is surjective we can choose $\xi$ so that the right side of~\eqref{eq:obs-change} is minus the original obstruction. Set $b_{2,R'}=b^0_{2,R'}+\xi$. The obstruction vanishes, so there is a perfect lift $\widetilde P_{R'}$ on $Z_{R'}=\mathcal X_{1,R'}\times_{R'}\mathcal X_{2,R'}$.

\smallskip\noindent\emph{Step 4: the lift is an equivalence.} Apply the relative projection to the kernel component $\AA^\vee_{1,R'}\boxtimes_{R'}\AA_{2,R'}\subset\Perf(Z_{R'})$, and call the result $P_{R'}$. By base change, its reduction is $P_R$, because $P_R$ already lies in that component. Let $F_{R'}$ be the induced functor and $P^R_{R'}$ the projected kernel of its right adjoint; it is perfect, since the ambient adjoint is given by relative duality for the smooth projective families, followed by the perfect projections. Let $Q_{i,R'}$ be the identity kernel of $\AA_{i,R'}$, realized as the right projection kernel. With $\star$ denoting convolution, the unit and counit of the adjunction are
\begin{equation}\label{eq:unit-counit}
 Q_{1,R'}\longrightarrow P^R_{R'}\star P_{R'},\qquad P_{R'}\star P^R_{R'}\longrightarrow Q_{2,R'}.
\end{equation}
Their cones are perfect, and their closed fibres vanish because $P$ represents an equivalence. The underlying topological spaces of these products are their closed fibres, and at each point a perfect complex with zero derived reduction to the residue field is acyclic; by derived Nakayama both cones vanish. So both maps in~\eqref{eq:unit-counit} are isomorphisms and $F_{R'}$ is an equivalence.

Interchanging the factors and using surjectivity of $\kappa_1$ gives the other assertion. Finally, every surjection of local Artinian $\C$-algebras with nilpotent kernel factors into small extensions, and iterating the construction proves the last statement.
\end{proof}

\begin{remark}\label{rem:one-side}
For a prescribed deformation of the first factor, the proof uses only the surjectivity of $\kappa_2$. Once both geometric lifts are chosen and the obstruction vanishes, the isomorphism classes of perfect lifts of the kernel form a torsor under $\Ext^1(P,P)\otimes J$.
\end{remark}

\subsection{Algebraic correspondences and specialization}\label{sec:special}

Theorem~\ref{thm:lift} is an infinitesimal statement. We now realize it by smooth morphisms from an algebraic variety carrying an actual family of equivalence kernels, and use this to spread properties of general fibres to all fibres. The vanishing of negative Hochschild cohomology enters here, through the algebraicity of the stack of complexes.

\begin{proposition}\label{prop:correspondence}
Assume the relative hypotheses of Subsection~\ref{sec:kernels}, and let $P$ be a projected equivalence kernel at $(b_1,b_2)$. Suppose that
\begin{equation}\label{eq:hyp-corr}
 \kappa_i\ \text{is surjective and}\ \HH^m(\AA_{i,0})=0\ \ (m<0),\qquad i=1,2.
\end{equation}
Then there are a smooth irreducible complex variety $T$, a point $t\in T(\C)$, smooth morphisms $p_i:T\to B_i$ with $p_i(t)=b_i$, and a perfect kernel on $(\mathcal X_1\times_{B_1}T)\times_T(\mathcal X_2\times_{B_2}T)$ inducing an equivalence of the relative components, whose fibre at $t$ is isomorphic to $P$.
\end{proposition}

\begin{proof}
\emph{Step 1: the stack of complexes.} Let $B=B_1\times B_2$ and $\mathcal Z=(\mathcal X_1\times B_2)\times_B(B_1\times\mathcal X_2)$, a smooth projective family over $B$. By \cite[Theorem~4.2.1]{Lie06}, the relatively perfect universally gluable complexes on $\mathcal Z/B$ form an algebraic stack $\mathcal P$, locally of finite presentation over $B$. Universal gluability means that negative relative self-Ext groups vanish after every base change; for relatively perfect complexes in this proper flat family it can be checked on geometric fibres \cite[Proposition~2.1.9]{Lie06}. By Lemma~\ref{lem:ext-hh} and~\eqref{eq:hyp-corr}, $P$ defines a point $e\in\mathcal P(\C)$.

\smallskip\noindent\emph{Step 2: the open substack of equivalence kernels.} Take a scheme chart $S\to\mathcal P$ with corresponding complex $C$ on $\mathcal Z_S$; since $\mathcal Z_S\to S$ is smooth, $C$ is perfect. Let $\Pi_S$ be the right projection to the relative kernel component. The cone of the counit $\Pi_S(C)\to C$ is perfect and commutes with base change, and its fibre at $s$ vanishes exactly when $C_s$ lies in the required component. For a perfect complex $D$ on a proper family $W\to S$, the locus $\{s\in S:D_s=0\}$ is open: the cohomological support of $D$ is closed, its image in $S$ is closed by properness, and derived Nakayama identifies the complement of the image with this locus. So the projected kernels form an open locus. On it, form the relative right adjoint kernel and the maps~\eqref{eq:unit-counit}; the same argument gives an open locus where both are isomorphisms. These loci are invariant under isomorphism and compatible with pullback, so they descend to an open substack $\mathcal E\subset\mathcal P$ of equivalence kernels containing $e$.

\smallskip\noindent\emph{Step 3: smoothness.} Choose a smooth scheme atlas $u:T\to\mathcal E$, of finite type near a point $t\in T(\C)$ above $e$, and let $p_i:T\to B_i$ be the composite maps. To see that $p_i$ is smooth at $t$, consider a small extension $R'\to R$ and compatible morphisms $\Spec R\to T$, $\Spec R'\to B_i$ with closed point $t$. The induced point of $\mathcal E$ consists of deformations of the two fibres and a projected equivalence kernel. Theorem~\ref{thm:lift} lifts this point to $\Spec R'$ while keeping the prescribed $B_i$-valued point. Since $u$ is smooth, the lifted point lifts to $T$ \'etale locally, and an \'etale neighbourhood of the closed point of $\Spec R'$ has a section, so the lift exists over $R'$. By the infinitesimal criterion \cite[Tag~02HX]{Stacks}, $p_i$ is smooth at $t$. Shrinking $T$ around $t$, both $p_i$ are smooth, and so is $T$, since the $B_i$ are smooth; we replace $T$ by its irreducible component through $t$, which is open. The object defining $T\to\mathcal E$ is the required family of kernels.
\end{proof}

\begin{theorem}\label{thm:specialize}
Assume the hypotheses of Proposition~\ref{prop:correspondence}, with $B_1$ and $B_2$ irreducible. Let $c_i:B_i\to\mathcal M$ be morphisms to a separated algebraic space such that equality of values at geometric points is equivalent to polarized isomorphism of the fibres. Suppose that there are dense open subsets $U_i\subset B_i$ with
\begin{equation}\label{eq:generic-hyp}
 u_i\in U_i(\C),\ \AA_{1,u_1}\simeq\AA_{2,u_2}\ \Longrightarrow\ c_1(u_1)=c_2(u_2).
\end{equation}
Then $c_1(b_1)=c_2(b_2)$; in particular, the original fibres are isomorphic as polarized varieties.
\end{theorem}

\begin{proof}
Let $T$, $t$ and $p_i$ be as in Proposition~\ref{prop:correspondence}. Each $p_i$ is smooth, hence open, and its image is a nonempty open subset of the irreducible variety $B_i$; so $p_i$ is dominant, and $T^\circ=p_1^{-1}(U_1)\cap p_2^{-1}(U_2)$ is dense and open in $T$. The family of kernels gives equivalences at all complex points of $T$, so by~\eqref{eq:generic-hyp} the maps $c_1p_1$ and $c_2p_2$ agree on the complex points of $T^\circ$. Since $\mathcal M$ is separated, $D=(c_1p_1,c_2p_2)^{-1}(\Delta_{\mathcal M})$ is closed in $T$; it contains $T^\circ$, hence $t$.
\end{proof}

\begin{remark}\label{rem:one-sided-generic}
If the generic statement only requires $u_1\in U_1$ and allows arbitrary $u_2\in B_2$, the proof works with $p_1^{-1}(U_1)$ in place of $T^\circ$. This is the form used for quartics below.
\end{remark}

\begin{remark}\label{rem:very-general}
Theorem~\ref{thm:specialize} also holds if each $B_i\setminus U_i$ is a countable union of proper closed subsets. Pull these subsets back to $T$; dominance of $p_i$ makes every pullback proper. A proper closed algebraic subset of the smooth irreducible variety $T$ has empty interior in the analytic topology, so by the Baire category theorem, applied on an analytic coordinate neighbourhood inside any nonempty Zariski open subset, the locus where both fibres lie in their $U_i$ is Zariski dense. The closed-diagonal argument then applies. Here the hypothesis must hold for all pairs $(u_1,u_2)$ with $u_i$ in the individual very-general loci.
\end{remark}

\subsubsection{Specialization of a uniform bound}
The same argument transports a uniform bound on the number of partners. Several kernels are deformed at once, and separatedness keeps their target fibres pairwise nonisomorphic.

\begin{theorem}\label{thm:uniform-special}
Let $\mathcal X\to B$ be a smooth projective polarized family over a smooth irreducible complex variety, with relative admissible components and perfect projection kernels compatible with base change. Suppose that the categorical Kodaira--Spencer map is surjective and negative Hochschild cohomology vanishes at every fibre. Let $c_B:B\to\mathcal M$ be a morphism to a separated algebraic space which distinguishes polarized isomorphism classes of fibres. Suppose that an integer $N$ bounds the size of every set of pairwise nonisomorphic fibres over a nonempty open subset $U\subset B$ whose components are mutually equivalent. Then the same bound holds over $B$. The same is true if $B\setminus U$ is a countable union of proper closed subsets, provided $N$ does not depend on the fibres in $U$.
\end{theorem}

\begin{proof}
Suppose that $\mathcal X_b$ has $N+1$ pairwise nonisomorphic partners $\mathcal X_{b_1},\dots,\mathcal X_{b_{N+1}}$. Choose an equivalence from the component of $\mathcal X_b$ to each of theirs and take its projected kernel. Proposition~\ref{prop:correspondence} gives smooth correspondences $T_i$ through points $t_i$, with smooth maps $p_i,q_i:T_i\to B$, $p_i(t_i)=b$, $q_i(t_i)=b_i$, and families of equivalence kernels.

Form $\widetilde T=T_1\times_B\dots\times_BT_{N+1}$ using the maps $p_i$; it contains $\widetilde t=(t_1,\dots,t_{N+1})$. The morphism $\widetilde T\to B$ and the projections $\widetilde T\to T_i$ are smooth, so $\widetilde T$ is smooth; let $T$ be its irreducible component through $\widetilde t$, which is open, and $t$ the distinguished point. The maps $r_i:T\to T_i\xrightarrow{q_i}B$ are smooth with nonempty open image, hence dominant. For $i\neq j$ the locus $E_{ij}=(c_Br_i,c_Br_j)^{-1}(\Delta_{\mathcal M})$ is closed, since $\mathcal M$ is separated, and does not contain $t$, since $\mathcal X_{b_i}\not\simeq\mathcal X_{b_j}$. So $T^{\mathrm{dist}}=T\setminus\bigcup_{i<j}E_{ij}$ is a nonempty open subset on which the target fibres are pairwise nonisomorphic.

If $U$ is open, each $r_i^{-1}(U)$ is dense and open, and their intersection with $T^{\mathrm{dist}}$ is nonempty. At a complex point of this intersection, the $N+1$ target fibres lie in $U$, are pairwise nonisomorphic, and have mutually equivalent components, since the pulled-back kernels identify each of them with the component of the common source fibre. This contradicts the bound on $U$. If $B\setminus U$ is a countable union of proper closed subsets, their preimages under the dominant maps $r_i$ are proper closed subsets of $T$, and the Baire argument of Remark~\ref{rem:very-general}, applied inside $T^{\mathrm{dist}}$, gives the same contradiction.
\end{proof}

\begin{remark}\label{rem:uniform-needed}
The argument needs one integer $N$ that bounds all partner sets over $U$. Proposition~\ref{prop:uniform} provides such an integer, independent of the target algebra.
\end{remark}

\begin{proof}[Proof of Theorem~\ref{thm:deformation}]
Part~(i) is Proposition~\ref{prop:correspondence}, which uses Theorem~\ref{thm:lift}. Part~(ii) is Theorem~\ref{thm:uniform-special}; its case $N=1$ is the Torelli statement.
\end{proof}

\section{Deformations and Hochschild cohomology of Kuznetsov components}\label{sec:hochschild}

In this section we verify the hypotheses of Theorem~\ref{thm:deformation} for the universal family of smooth hypersurfaces, and identify the intrinsic algebra used in the Calabi--Yau case. Subsection~\ref{sec:def-hyp} treats hypersurfaces of degree at least four, the next subsection the universal family, and Subsection~\ref{sec:cy} the Calabi--Yau case, including cubics.

\subsection{Deformations of hypersurfaces of degree at least four}\label{sec:def-hyp}

We now verify the hypotheses of Theorem~\ref{thm:deformation} for hypersurfaces of degree at least four: the geometric deformations surject onto $\HH^2(\Ku(X))$, and negative Hochschild cohomology vanishes. Let $X=V(f)\subset\PP^{M-1}$ be a smooth hypersurface of degree $d$, with
\begin{equation}\label{eq:range4}
 4\leq d\leq M-1.
\end{equation}
Put $n=M-2$ and $I=M-d$, so that $n\geq3$, $I>0$ and $\omega_X^{-1}\simeq\OO_X(I)$. We use the Kuznetsov component~\eqref{eq:ku} with its induced dg enhancement, and write
\[
 S=\C[x_1,\dots,x_M],\qquad J_f=(\partial_1f,\dots,\partial_Mf),\qquad R_f=S/J_f.
\]

\begin{theorem}\label{thm:deg4}
For every smooth hypersurface satisfying~\eqref{eq:range4}, the map $\eta_X$ is an isomorphism and $\HH^m(\Ku(X))=0$ for $m<0$. Moreover,
\begin{equation}\label{eq:dim-hh2}
 \dim\HH^2(\Ku(X))=h^1(X,T_X)=\binom{M+d-1}{d}-M^2.
\end{equation}
\end{theorem}

The proof shows first that $\eta_X$ is injective, using Kuznetsov's normal Hochschild cohomology, and then that source and target have the same dimension.

\begin{lemma}\label{lem:pseudoheight}
The exceptional collection in~\eqref{eq:ku} has pseudoheight $d-1$. Consequently, the restriction map $\HH^2(X)\to\HH^2(\Ku(X))$ is injective, and $\HH^m(\Ku(X))=0$ for $m<0$.
\end{lemma}

\begin{proof}
Set $E_a=\OO_X(a)$ for $0\leq a\leq I-1$, and let $S_X$ be the Serre functor of $\Perf(X)$. For $a<b$, Kodaira vanishing gives $\Ext^q(E_a,E_b)=H^q(X,\OO_X(b-a))=0$ for $q>0$, whereas $\Hom(E_a,E_b)\neq0$; the same vanishing and anti-ampleness give semiorthogonality when $-I<b-a<0$. So $(E_0,\dots,E_{I-1})$ is exceptional. For a chain $0\leq a_0<\dots<a_p\leq I-1$, the closing term in the definition of pseudoheight is
\begin{equation}\label{eq:closing}
 \RHom(E_{a_p},S_X^{-1}E_{a_0})\simeq R\Gamma\bigl(X,\OO_X(I+a_0-a_p)\bigr)[-n].
\end{equation}
The twist is positive, so the right side has cohomology only in degree $n$. All forward terms have their first nonzero Ext group in degree zero. By the definition of pseudoheight \cite[Definition~4.4]{Kuz15}, this chain contributes $n-p$. Since $p\leq I-1$, every contribution is at least $n-I+1=d-1$, and the full chain $0<1<\dots<I-1$ attains this value, including the length-zero chain when $I=1$. So the pseudoheight is $d-1$.

The normal Hochschild triangle \cite[Theorem~3.3]{Kuz15} and the pseudoheight bound \cite[Lemma~4.5]{Kuz15} give $\NHH^q(\langle E_0,\dots,E_{I-1}\rangle,\Perf(X))=0$ for $q<d-1$. Since $d-1\geq3$, restriction is an isomorphism in degrees at most $d-3$ and injective in degree $d-2$ \cite[Corollary~4.6]{Kuz15}, in particular injective in degree two. The Hochschild--Kostant--Rosenberg theorem gives $\HH^m(X)=0$ for $m<0$, and restriction transfers this vanishing to $\Ku(X)$.
\end{proof}

\begin{proof}[Proof of Theorem~\ref{thm:deg4}]
\emph{Injectivity.} The Hochschild--Kostant--Rosenberg decomposition \cite{Swa96} is
\[
 \HH^2(X)\simeq H^0(X,\textstyle\bigwedge^2T_X)\oplus H^1(X,T_X)\oplus H^2(X,\OO_X),
\]
so the geometric summand is included injectively, and Lemma~\ref{lem:pseudoheight} shows that $\eta_X$ is injective.

\smallskip\noindent\emph{The source.} The restriction sequence $0\to\OO_{\PP^{M-1}}(t-d)\to\OO_{\PP^{M-1}}(t)\to\OO_X(t)\to0$ gives $H^q(X,\OO_X(t))=0$ for $0<q<n$ and all $t$. The restricted Euler sequence $0\to\OO_X\to\OO_X(1)^{\oplus M}\to T_{\PP^{M-1}}|_X\to0$ then gives $H^1(X,T_{\PP^{M-1}}|_X)=0$, and the normal sequence yields
\begin{equation}\label{eq:normal}
 H^1(X,T_X)\simeq\coker\bigl(H^0(X,T_{\PP^{M-1}}|_X)\to H^0(X,\OO_X(d))\bigr).
\end{equation}
After choosing $f$, the target is $S_d/\C f$, and an ambient vector field given by a matrix $A=(a_{ij})$ maps to $D_Af=\sum_{i,j}a_{ij}x_i\partial_jf$ modulo $\C f$. By Euler's identity $f\in(J_f)_d$, so
\begin{equation}\label{eq:h1-jac}
 H^1(X,T_X)\simeq(R_f)_d.
\end{equation}
Since $X$ is smooth, the affine critical locus of $f$ is the origin, and the partial derivatives form a regular sequence. Their first Koszul syzygies have degree $2d-2>d$, so the $M^2$ products $x_i\partial_jf$ are linearly independent and
\begin{equation}\label{eq:dim-jac}
 \dim(R_f)_d=\dim S_d-M^2=\binom{M+d-1}{d}-M^2.
\end{equation}

\smallskip\noindent\emph{The target.} We use the dg equivalence of $\Ku(X)$ with graded matrix factorizations \cite[Theorem~2.1]{LRZ26}, obtained from \cite[Theorem~6.13]{BFK14}. Denote the diagonal bimodule by $\Delta$ and the internal degree shift by $(1)$, with
\begin{equation}\label{eq:shift}
 (d)\simeq[2].
\end{equation}
For odd $M$, \cite[Proposition~3.10]{LRZ26} gives $\Hom(\Delta,\Delta(t))\simeq(R_f)_t$. For even $M$ it gives
\begin{equation}\label{eq:sectors}
 \Hom(\Delta,\Delta(t))\simeq
 \begin{cases}(R_f)_t,&t\neq\frac{(d-2)M}{2},\\[2pt](R_f)_t\oplus\C^{d-1},&t=\frac{(d-2)M}{2}.\end{cases}
\end{equation}
When $M$ is even, $\frac{(d-2)M}{2}\geq M>d$, because $d\geq4$ and $M\geq d+1$, so the extra summand does not occur in internal degree $d$. In both cases
\begin{equation}\label{eq:hh2-jac}
 \HH^2(\Ku(X))\simeq\Hom(\Delta,\Delta[2])\simeq\Hom(\Delta,\Delta(d))\simeq(R_f)_d.
\end{equation}
By~\eqref{eq:h1-jac} and~\eqref{eq:hh2-jac}, the injective map $\eta_X$ has source and target of the same finite dimension, so it is an isomorphism, and~\eqref{eq:dim-hh2} follows from~\eqref{eq:dim-jac}.
\end{proof}

\begin{proposition}\label{prop:quartic-hh}
Let $X=V(f)\subset\PP^{4k+1}$ be a smooth quartic, $k\geq1$. For every integer $a$,
\begin{equation}\label{eq:quartic-hh}
 \HH^{2a}(\Ku(X))\simeq(R_f)_{4a},\qquad\HH^{2a+1}(\Ku(X))=0.
\end{equation}
In particular, $H^1(X,T_X)\to\HH^2(\Ku(X))$ is an isomorphism of dimension $\binom{4k+5}{4}-(4k+2)^2$.
\end{proposition}

\begin{proof}
Here $M=4k+2$ is even and $[2a]\simeq(4a)$ by~\eqref{eq:shift}. The exceptional internal degree in~\eqref{eq:sectors} is $M\equiv2\pmod4$, so it never equals $4a$; this gives the even degrees. For even $M$, \cite[Proposition~3.10(1)]{LRZ26} also gives $\Hom(\Delta,\Delta(t)[1])=0$ for all $t$, which gives the odd degrees. The last assertion is Theorem~\ref{thm:deg4}.
\end{proof}

\subsection{The universal family of hypersurfaces}
Let $B_{d,M}\subset\PP(S_d)$ be the open set of smooth equations and $\pi:\mathcal X\to B_{d,M}$ the universal hypersurface. The base is smooth and irreducible. We polarize the family by $\OO_{\mathcal X}(1)$ and use the relative Kuznetsov component defined by $\OO_{\mathcal X},\dots,\OO_{\mathcal X}(I-1)$.

\begin{lemma}\label{lem:equation-family}
At every $b=[f]\in B_{d,M}$, the map $\kappa_b:T_bB_{d,M}\to\HH^2(\Ku(\mathcal X_b))$ is surjective. The equation family is locally versal for polarized deformations, and the relative Kuznetsov component has perfect projection kernels compatible with base change.
\end{lemma}

\begin{proof}
After choosing $f$, the tangent space of the equation space is $S_d/\C f$, and its Kodaira--Spencer map is the quotient onto $S_d/(J_f)_d$ in~\eqref{eq:normal}. It is surjective, and so is its composite with the isomorphism $\eta_{\mathcal X_b}$ of Theorem~\ref{thm:deg4}.

For local versality: since $H^1(X,\OO_X)=H^2(X,\OO_X)=0$, the line bundle $\OO_X(1)$ extends uniquely through a small extension, and $H^1(X,\OO_X(1))=0$ lets us lift a basis of its sections. These sections give a closed immersion into projective space over the thickening, since the relative complete linear system stays very ample. The image is a flat effective Cartier divisor whose ideal reduces to $\OO_{\PP^{M-1}}(-d)$; since line bundles on projective space have no infinitesimal deformations, the ideal is $\OO(-d)$ and the image is defined by a homogeneous equation of degree $d$. Lifting bases compatibly at each step gives local versality.

The relative line bundles restrict to the exceptional collection on every fibre; the vanishings in Lemma~\ref{lem:pseudoheight} and cohomology and base change show that they form a relative exceptional collection. The relative left mutation through $\OO_{\mathcal X}(j)$ has kernel
\[
 \Cone\bigl(p_1^*\OO_{\mathcal X}(-j)\otimes p_2^*\OO_{\mathcal X}(j)\to\OO_{\Delta_{\mathcal X/B_{d,M}}}\bigr)
\]
on $\mathcal X\times_{B_{d,M}}\mathcal X$, which is perfect and commutes with derived base change. Composing these mutations gives the left projection to the relative component, and its right adjoint gives the right projection; relative duality and smooth properness show that the adjoint kernels are perfect and commute with base change.
\end{proof}

\begin{corollary}\label{cor:generic-to-all}
Fix $d,M$ satisfying~\eqref{eq:range4}. Suppose that there is a nonempty Zariski open subset $U\subset B_{d,M}$ such that, for $u,v\in U(\C)$, a Fourier--Mukai equivalence $\Ku(\mathcal X_u)\simeq\Ku(\mathcal X_v)$ implies a projective isomorphism $\mathcal X_u\simeq\mathcal X_v$. Then the same holds for every pair of smooth degree-$d$ hypersurfaces in $\PP^{M-1}$. The same conclusion holds if $U$ is an individual very-general locus.
\end{corollary}

\begin{proof}
Theorem~\ref{thm:deg4} and Lemma~\ref{lem:equation-family} verify the deformation and kernel hypotheses at every pair. By \cite[Theorem~1.7 and Corollary~1.8]{Ben13}, the moduli space of smooth degree-$d$ hypersurfaces polarized by $\OO(1)$ is a separated algebraic space that distinguishes polarized isomorphism classes; the quadric exception is excluded by $d\geq4$. A polarized isomorphism is induced by a projective linear isomorphism of the embeddings by $|\OO(1)|$. Theorem~\ref{thm:specialize} and Remark~\ref{rem:very-general}, applied to two copies of the universal family, give the assertion.
\end{proof}

\subsection{The Calabi--Yau case}\label{sec:cy}

In this subsection $d$ divides $M$, so the Kuznetsov component is a Calabi--Yau category. We identify an intrinsic algebra, the $d$-th Veronese of the Jacobian ring, and verify the deformation hypotheses for cubics, which are not covered by Subsection~\ref{sec:def-hyp}. We keep the notation~\eqref{eq:cy-range} and set
\[
 n=M-2,\qquad I=M-d,\qquad c=m(d-2),\qquad R_f=\C[x_1,\dots,x_M]/(\partial_1f,\dots,\partial_Mf).
\]
The partial derivatives form a regular sequence, so $R_f$ is an Artinian Gorenstein algebra with Hilbert series
\begin{equation}\label{eq:hilbert}
 \sum_{j\geq0}\dim(R_f)_jt^j=(1+t+\dots+t^{d-2})^M
\end{equation}
and socle degree $M(d-2)=cd$.

\subsubsection{Serre and rotation functors}
Let $\rho$ be the rotation functor of $\Ku(X)$, obtained by tensoring with $\OO_X(1)$ and mutating back into the component. The matrix-factorization description gives
\begin{equation}\label{eq:rotation}
 \rho^d\simeq[2],\qquad S_{\Ku(X)}\simeq\rho^{-M}[M];
\end{equation}
see \cite[Theorem~2.1 and Proposition~2.2]{LRZ26}, \cite{Orl09} and \cite{Kuz19}. Since $M=md$,
\begin{equation}\label{eq:cy-serre}
 S_{\Ku(X)}\simeq[c],\qquad c=M-2m=m(d-2),
\end{equation}
so $\Ku(X)$ is a Calabi--Yau category of dimension $c$. For cubics:
\begin{center}
\renewcommand{\arraystretch}{1.1}
\begin{tabular}{cccc}
\toprule
$\dim X$ & $M$ & $m$ & Calabi--Yau dimension of $\Ku(X)$\\
\midrule
4 & 6 & 2 & 2\\
7 & 9 & 3 & 3\\
10 & 12 & 4 & 4\\
13 & 15 & 5 & 5\\
\bottomrule
\end{tabular}
\end{center}
The first row, cubic fourfolds, is excluded from Theorem~\ref{thm:finite}.

For a smooth proper dg category $\AA$ with Serre functor $S_\AA$, the \emph{Hochschild--Serre groups} are $\HS^{p,q}(\AA)=\Hom(\Id_\AA,S_\AA^q[p])$. When $S\simeq[c]$, they reduce to Hochschild cohomology,
\begin{equation}\label{eq:hs-cy}
 \HS^{p,q}(\AA)\simeq\HH^{p+cq}(\AA),
\end{equation}
so the Hochschild--Serre algebra does not single out a $d$-th root of $[2]$.

\subsubsection{The Jacobian Veronese inside Hochschild cohomology}

\begin{proposition}\label{prop:veronese-intrinsic}
For every integer $r$,
\begin{equation}\label{eq:cy-hh}
 \HH^r(\Ku(X))\simeq\begin{cases}(R_f)_{dr/2},&r\ \text{even},\\0,&r\ \text{odd},\end{cases}\ \oplus\ \begin{cases}\C^{d-1},&r=c,\\0,&r\neq c.\end{cases}
\end{equation}
The unital subalgebra generated by $\HH^2(\Ku(X))$, placed in degree one, is canonically identified with
\begin{equation}\label{eq:BX}
 B_X=R_f^{(d)}=\bigoplus_{a\geq0}(R_f)_{ad}.
\end{equation}
It is a standard graded Gorenstein algebra with socle degree $c$. Consequently, an enhanced equivalence $\Ku(X)\simeq\Ku(Y)$ induces a graded algebra isomorphism $B_X\simeq B_Y$.
\end{proposition}

\begin{proof}
Use the formula for the diagonal matrix factorization \cite[Theorem~5.39]{BFK14}, in the form of \cite[Theorem~3.5, Remarks~3.6--3.7, Propositions~3.10--3.11]{LRZ26}. The identity sector contributes $(R_f)_{dr/2}$ in even cohomological degree. Every nonidentity element of the scalar group $\mu_d$ fixes only the origin; its conormal rank is $M$ and its determinant character is $-M$. Its invariant contribution therefore requires $r\equiv M\pmod2$ and $M+\frac d2(r-M)=0$, which, as $d\mid M$, gives $r=M-2M/d=c$. There are $d-1$ nonidentity elements, which proves~\eqref{eq:cy-hh}.

In the range~\eqref{eq:cy-range}, $c=2$ only for $(d,m)=(3,2)$, so all of $\HH^2$ is the identity-sector space $(R_f)_d$. Products of identity-sector classes are Jacobian products. The $d$-th Veronese of a quotient of a polynomial ring is generated by its first piece, since every monomial of degree $ad$ splits into $a$ monomials of degree $d$; this proves the multiplicative assertion. The Gorenstein pairing of $R_f$ restricts to perfect pairings $(B_X)_a\otimes(B_X)_{c-a}\to(B_X)_c$. Finally, enhanced equivalences preserve the cup product and its grading, hence the subalgebra generated by $\HH^2$.
\end{proof}

\begin{remark}\label{rem:extra-sectors}
If $c$ is even, the additional sectors lie in an even Hochschild degree, but they are not in the subalgebra generated by $\HH^2$, since products of identity-sector classes stay in the identity sector. For a cubic sevenfold,
\[
 \bigl(\dim\HH^r(\Ku(X))\bigr)_{0\leq r\leq6}=(1,0,84,2,84,0,1).
\]
\end{remark}

\subsubsection{The deformation map for Calabi--Yau cubics}

\begin{proposition}\label{prop:cubic-def}
Let $m\geq3$ and let $X\subset\PP^{3m-1}$ be a smooth cubic. The map
\[
 \eta_X:H^1(X,T_X)\longrightarrow\HH^2(X)\longrightarrow\HH^2(\Ku(X))
\]
is an isomorphism, and $\HH^r(\Ku(X))=0$ for $r<0$ and for $r=1$.
\end{proposition}

\begin{proof}
Write $M=3m$ and $n=M-2$. The Euler and normal sequences give $H^1(X,T_X)\simeq(R_f)_3$, of dimension $\binom{M+2}{3}-M^2$. The exceptional collection is $\OO_X,\dots,\OO_X(M-4)$, and the primitive Hodge pieces
\begin{equation}\label{eq:cubic-hodge}
 H^{2m-1,m-1}_{\prim}(X)\simeq(R_f)_0,\qquad H^{2m-2,m}_{\prim}(X)\simeq(R_f)_3
\end{equation}
show that the Hodge diamond is not diagonal, so \cite[Lemma~5.2]{Pir24} applies to this collection and shows that $\eta_X$ is injective. In this case the argument is direct. With the convention $\HH_j(X)=\bigoplus_{p-q=j}H^q(X,\Omega^p_X)$, the two spaces in~\eqref{eq:cubic-hodge} lie in degrees $m$ and $m-2$, both nonzero, so the inclusion of the component identifies these parts of Hochschild homology with the corresponding parts for $X$. Under the residue description, a deformation $\xi\in(R_f)_3$ acts by $(R_f)_0\to(R_f)_3$, $1\mapsto\xi$. The compatibility of this action with restriction, including the Todd-corrected Hochschild--Kostant--Rosenberg identification, is established in the proof of \cite[Lemma~5.2]{Pir24}; the same identification of the composite is made in the proof of \cite[Proposition~3.5]{HR19}.

The matrix-factorization calculation gives
\begin{equation}\label{eq:cubic-hh}
 \HH^r(\Ku(X))\simeq\begin{cases}(R_f)_{3r/2},&r\ \text{even},\\0,&r\ \text{odd},\end{cases}\ \oplus\ \begin{cases}\C^2,&r=m,\\0,&r\neq m.\end{cases}
\end{equation}
The identity sector is the indicated Jacobian piece. For each nonidentity scalar in $\mu_3$ the fixed locus is the origin, and the invariant degree condition is $M+\frac32(r-M)=0$, that is, $r=m$, where the parity condition holds automatically; this is the case $d=3$ of \cite[Theorem~3.5 and Proposition~3.10]{LRZ26}. Since $m\geq3$, $\dim\HH^2(\Ku(X))=\dim(R_f)_3$, so the injective map $\eta_X$ is an isomorphism. Formula~\eqref{eq:cubic-hh} also gives the vanishings.
\end{proof}

\begin{remark}\label{rem:cubic-fourfold}
For a cubic fourfold ($m=2$), the same calculation gives $\HH^2(\Ku(X))\simeq(R_f)_3\oplus\C^2$: the geometric deformation space has dimension $20$ and the categorical one has dimension $22$; see also \cite[Remark~3.8]{HR19}. The map $\kappa$ is not surjective, and Theorem~\ref{thm:deformation} does not apply.
\end{remark}

\begin{proposition}\label{prop:cy-hyp}
For every hypersurface satisfying~\eqref{eq:cy-range}, the map $\eta_X$ is an isomorphism, and $\HH^r(\Ku(X))=0$ for $r<0$ and $r=1$. The universal family of smooth equations has surjective maps $\kappa$ and perfect relative left and right projection kernels compatible with base change.
\end{proposition}

\begin{proof}
For $d\geq4$ the assertion about $\eta_X$ is Theorem~\ref{thm:deg4}; for $d=3$ the condition $(d,m)\neq(3,2)$ gives $m\geq3$, and Proposition~\ref{prop:cubic-def} applies. The vanishings follow from~\eqref{eq:cy-hh}, since $c\geq3$. On the space $B_{d,M}\subset\PP(S_d)$ of smooth equations, the Euler and normal sequences identify the Kodaira--Spencer map with the surjection $S_d/\C f\to(R_f)_d=H^1(X_f,T_{X_f})$, so its composite with $\eta_{X_f}$ is surjective. The construction of the relative projection kernels in Lemma~\ref{lem:equation-family} applies verbatim, including for cubics.
\end{proof}

\section{Quartic hypersurfaces}\label{sec:quartics}

In this section we prove Theorem~\ref{thm:quartic}. Subsection~\ref{sec:even-jac} recovers a quartic from its even Jacobian algebra under a genericity condition, Subsection~\ref{sec:generic-quartic} shows that this condition holds for a general quartic, and Subsection~\ref{sec:quartic-torelli} recovers the even Jacobian algebra from the Kuznetsov component and applies Theorem~\ref{thm:deformation}.

\subsection{Reconstruction from the even Jacobian algebra}\label{sec:even-jac}

For quartics in $\PP^{4k+1}$, the category gives access only to the even part of the Jacobian ring (Subsection~\ref{sec:quartic-torelli}). This subsection shows how to recover the polynomial structure from it. The idea is that the relations among products of quadrics coming from the polynomial ring itself are exactly the relations of rank at most three, and that, for a quartic satisfying a genericity condition, the Jacobian ideal contributes no further relations of rank at most three.

Let $V$ be a complex vector space of dimension $M\geq3$, $S=\Sym V$, and $S_j=\Sym^jV$; so $V$ is the space of linear forms on $\PP(V^*)$. For $F\in S_4$ let $J_F=(\partial_vF:v\in V^*)$ and $R_F=S/J_F$. The \emph{even Jacobian algebra} is
\begin{equation}\label{eq:even-alg}
 R^{(2)}_F=\bigoplus_{a\geq0}(R_F)_{2a},\qquad(R^{(2)}_F)_a=(R_F)_{2a}.
\end{equation}
Put $W=S_2$ and $q=\dim W=\binom{M+1}{2}$. Multiplication gives
\begin{equation}\label{eq:mu}
 \mu:\Sym^2W\longrightarrow S_4,\qquad K_0=\ker\mu.
\end{equation}
We write $\odot$ for multiplication in $\Sym(W)$. An element of $\Sym^2W$ is a quadratic form on $W^*$, so its rank is defined intrinsically in terms of $W$.

\begin{lemma}\label{lem:veronese}
The homogeneous ideal of the quadratic Veronese embedding $\PP(V^*)\hookrightarrow\PP(W^*)$ is generated by $K_0$, and
\[
 K_0=\Span_\C\{Q\in K_0:\rank Q\leq3\}.
\]
\end{lemma}

\begin{proof}
This is the quadratic case of \cite[Theorem~3.1]{HLMP21}; we give a direct proof over $\C$. For $u,v\in V$ set $Q(u,v)=u^2\odot v^2-(uv)\odot(uv)$, an element of $K_0$ of rank at most three, and let $L$ be the span of these elements. For $a,b,c,e\in V$ write
\[
 A=(ab)\odot(ce),\qquad B=(ac)\odot(be),\qquad C=(ae)\odot(bc).
\]
The coefficient of $s\,t\,r\,w$ in $Q(sa+tb,rc+we)$ is $2(2A-B-C)$, so $2A-B-C\in L$, because a polynomial with values in a fixed subspace has all its coefficients in that subspace. Interchanging $b$ and $c$ gives $2B-A-C\in L$, and subtracting gives $3(A-B)\in L$, so
\begin{equation}\label{eq:exchange}
 (ab)\odot(ce)-(ac)\odot(be)\in L.
\end{equation}
Choose a basis $x_1,\dots,x_M$ of $V$. The homomorphism $\Sym(W)\to\bigoplus_{a\geq0}S_{2a}$ sends monomials in the symbols $x_ix_j$ to monomials in the $x_i$, so its kernel is spanned by differences of monomials with the same image, that is, by pairs of pairings of the same multiset of $2a$ indices. Any two pairings are connected by successively replacing pairs $(i,j),(l,p)$ by $(i,l),(j,p)$: match the pair containing the first index in one replacement, and proceed by induction on the remaining indices. Each replacement is a multiple of a relation~\eqref{eq:exchange}. So these quadratic relations generate the kernel, their span is $K_0$, and since they lie in $L\subset K_0$ the equality follows.
\end{proof}

Put
\begin{equation}\label{eq:Q3}
 \mathcal Q_3=\{Q\in\Sym^2W:\rank Q\leq3\},\qquad\mathcal C_3=\mu(\mathcal Q_3)\subset S_4.
\end{equation}
The set $\mathcal C_3$ is constructible. Over $\C$ every element of $\mathcal C_3$ has the form $ab-c^2$ with $a,b,c\in S_2$. Choosing a three-dimensional support for a symmetric form takes $3(q-3)$ parameters and choosing the form on it takes six; the lower-rank strata are smaller, so
\begin{equation}\label{eq:dimQ3}
 \dim\mathcal Q_3=3q-3.
\end{equation}

\begin{definition}\label{def:avoid}
A quartic $F\in S_4$ has the \emph{avoidance property} if $(J_F)_4\cap\mathcal C_3=\{0\}$.
\end{definition}

\begin{proposition}\label{prop:avoid-recon}
Let $F,G\in S_4$ define smooth hypersurfaces in $\PP(V^*)$, and suppose that $F$ has the avoidance property. If $R^{(2)}_F\simeq R^{(2)}_G$ as graded $\C$-algebras, then $V(F)$ and $V(G)$ are projectively isomorphic.
\end{proposition}

\begin{proof}
\emph{Step 1: the Veronese relations.} The generators of a quartic Jacobian ideal are cubics, so $(R_F)_2=S_2=W$. Multiplication in degree two determines
\[
 K_F=\ker\bigl(\Sym^2(R_F)_2\to(R_F)_4\bigr)=\mu^{-1}\bigl((J_F)_4\bigr).
\]
For a subspace $K\subset\Sym^2W$ let $L_3(K)$ be the span of its elements of rank at most three. The avoidance property gives $K_F\cap\mathcal Q_3=K_0\cap\mathcal Q_3$, so by Lemma~\ref{lem:veronese}
\begin{equation}\label{eq:L3}
 L_3(K_F)=K_0.
\end{equation}
Let $\psi:R^{(2)}_F\to R^{(2)}_G$ be an isomorphism. Its degree-one part identifies $K_F$ with $K_G$ and preserves ranks of quadratic forms, hence identifies $L_3(K_F)$ with $L_3(K_G)$. On the $G$ side, the Veronese relation space is contained in $L_3(K_G)$ by Lemma~\ref{lem:veronese}, and the two have the same dimension by~\eqref{eq:L3}, so they are equal.

\smallskip\noindent\emph{Step 2: the linear change of coordinates.} So $\psi_1$ induces a graded isomorphism of the quadratic Veronese coordinate ring $\Sym(W)/(K_0)\simeq\bigoplus_aS_{2a}$ with itself. Taking $\Proj$ gives an automorphism of projective space preserving $\OO(2)$, which is projective linear; hence $\psi_1$ is, up to a nonzero scalar, the second symmetric power of a linear automorphism of $V$. Concretely, after removing that automorphism, the remaining linear map of $W^*$ has every line of the Veronese variety as an eigenline; the variety is irreducible and spans $\PP(W^*)$, so the map is scalar, and over $\C$ a rescaling of the automorphism of $V$ absorbs this scalar.

\smallskip\noindent\emph{Step 3: the Jacobian ideals.} After this change of coordinates $\psi_1$ is the identity on $S_2$. The even Jacobian algebras are generated in degree one, so their degree-two multiplication kernels agree, and since $\mu$ is surjective, $(J_F)_4=(J_G)_4$. The partial derivatives of a smooth quartic form a regular sequence of $M$ cubics, so $R_F$ is an Artinian complete intersection with socle in degree $2M$; in particular it has no socle in degree three, and
\begin{equation}\label{eq:J3}
 (J_F)_3=\{h\in S_3:S_1h\subset(J_F)_4\}.
\end{equation}
The same holds for $G$. Hence $J_F=J_G$, and Donagi's reconstruction theorem for homogeneous Jacobian ideals \cite[Proposition~1.1]{Don83} shows that $V(F)$ and $V(G)$ are projectively isomorphic.
\end{proof}

\subsection{Generic reconstruction of quartics}\label{sec:generic-quartic}

Let $B_M\subset\PP(S_4)$ be the open set of smooth quartics. In this subsection we prove that the avoidance property holds generically, which gives the following theorem.

\begin{theorem}\label{thm:generic-quartic}
Let $M=4k+2$ with $k\geq1$. There is a nonempty Zariski open subset $U_M\subset B_M$ such that, for $[F]\in U_M$ and any $[G]\in B_M$, a graded algebra isomorphism $R^{(2)}_F\simeq R^{(2)}_G$ implies that $V(F)$ and $V(G)$ are projectively isomorphic.
\end{theorem}

By Proposition~\ref{prop:avoid-recon} it suffices to show that the avoidance property holds on a nonempty open subset of $S_4$. The estimates use the description of $(J_F)_4$ as the image of the linear derivations, given in~\eqref{eq:DA}.

\subsubsection{The incidence estimate}
For $A\in\End(V)$, let $D_A$ be the degree-zero derivation of $S$ induced by $A$. In a basis,
\begin{equation}\label{eq:DA}
 D_AF=\sum_{i,j=1}^Ma_{ij}x_i\frac{\partial F}{\partial x_j},\qquad(J_F)_4=\{D_AF:A\in\End(V)\}.
\end{equation}
Write $H=\dim S_4=\binom{M+3}{4}$.

\begin{lemma}\label{lem:incidence}
Let $T\subset\PP(\End V)$ be a locally closed subvariety of dimension at most $p$ such that $\rank(D_A:S_4\to S_4)\geq r$ for $[A]\in T$. If $r>p+3q-3$, the set of quartics $F$ with $D_AF\in\mathcal C_3$ for some $[A]\in T$ has proper Zariski closure in $S_4$.
\end{lemma}

\begin{proof}
Cover $T$ by finitely many affine charts with a chosen representative of $[A]$, and on each consider the incidence variety $I=\{(A,F,Q):Q\in\mathcal Q_3,\ D_AF=\mu(Q)\}$. For fixed $A$ and $Q$, the possible $F$ form an empty set or a translate of $\ker D_A$, of dimension at most $H-r$. By~\eqref{eq:dimQ3}, $\dim I\leq p+(3q-3)+H-r<H$, so the image in $S_4$ has proper closure. Take the finite union over the charts.
\end{proof}

\begin{lemma}\label{lem:rank-bound}
For every nonzero $A\in\End(V)$, $\rank(D_A:S_4\to S_4)\geq\binom{M+2}{3}$.
\end{lemma}

\begin{proof}
A nonzero scalar acts on $S_4$ by a nonzero scalar. If $A$ is not scalar, some off-diagonal entry $a_{ij}$ is nonzero in a suitable basis. Conjugating by the diagonal matrix with weights $1,-1,0,\dots,0$ at $i,j$ and the other indices, dividing by $t^2$ and letting $t\to\infty$, the matrix tends to $a_{ij}E_{ij}$, the unique entry of weight difference two. Conjugation and scaling preserve the rank of $D_A$, and rank does not increase under specialization, so it suffices to compute the rank of $x_i\partial_{x_j}$. Differentiation $S_4\to S_3$ is surjective and multiplication by $x_i$ is injective, so this rank is $\dim S_3=\binom{M+2}{3}$.
\end{proof}

\begin{proposition}\label{prop:M13}
If $M\geq13$, a quartic in a nonempty Zariski open subset of $S_4$ has the avoidance property.
\end{proposition}

\begin{proof}
Apply Lemma~\ref{lem:incidence} with $p=M^2-1$ and $r=\binom{M+2}{3}$. The difference of the two sides of its inequality is
\[
 \binom{M+2}{3}-(M^2-1)-\Bigl(3\binom{M+1}{2}-3\Bigr)=\frac{M^3-12M^2-7M+24}{6},
\]
which equals $17$ at $M=13$ and increases for $M\geq13$. So the quartics with $D_AF\in\mathcal C_3$ for some nonzero $A$ lie in a proper closed subset, and outside it~\eqref{eq:DA} gives the avoidance property.
\end{proof}

For $M=4k+2$ this treats $k\geq3$. The cases $M=10$ and $M=6$ need finer estimates.

\subsubsection{Ten variables}

\begin{lemma}\label{lem:subspace}
Let $\dim V\geq4$. If $A$ is not a scalar plus an endomorphism of rank at most one, there is a two-dimensional subspace $U\subset V$ with $U\cap AU=0$ and $\dim AU=2$.
\end{lemma}

\begin{proof}
If the minimal polynomial of $A$ has degree at least three, choose $v$ with $v,Av,A^2v$ independent and $w$ outside their span. For $U_\epsilon=\langle v,Av+\epsilon w\rangle$,
\[
 v\wedge(Av+\epsilon w)\wedge Av\wedge(A^2v+\epsilon Aw)=\epsilon\,v\wedge w\wedge Av\wedge A^2v+O(\epsilon^2),
\]
whose linear coefficient is nonzero, so $U_\epsilon$ works for general $\epsilon$. If the minimal polynomial has degree two with distinct roots, $A$ has two eigenspaces, each of dimension at least two since otherwise $A$ would be a scalar plus a rank-one endomorphism; with $e_1,e_2$ independent in one and $f_1,f_2$ in the other, $U=\langle e_1+f_1,e_2+f_2\rangle$ works. If the minimal polynomial has a repeated root, $A=a\,\mathrm{id}+N$ with $N^2=0$ and $\rank N\geq2$, and $U$ spanned by the top vectors of two independent Jordan chains of length two works. The scalar case is excluded.
\end{proof}

\begin{lemma}\label{lem:rank348}
Let $M=10$, and suppose that $A$ is not a scalar plus an endomorphism of rank at most one. Then $\rank(D_A:S_4\to S_4)\geq348$.
\end{lemma}

\begin{proof}
Choose $U$ as in Lemma~\ref{lem:subspace} and decompose $V=AU\oplus W'\oplus U$. Conjugating with weights $1,0,-1$ on these summands and scaling by $t^{-2}$ degenerates $A$ to the endomorphism that is an isomorphism $U\to AU$ and zero on the other summands, a square-zero nilpotent of Jordan type $(2,2,1^6)$. Complete it to an $\mathfrak{sl}_2$-action on two copies of the standard representation and six trivial ones; the weights on $V$ are $1,1,-1,-1,0,0,0,0,0,0$. In a finite-dimensional $\mathfrak{sl}_2$-representation, each irreducible summand contributes one dimension to the kernel of the raising operator and one dimension to exactly one of the weight spaces of weight $0$ and $1$. So the kernel on $S_4$ has dimension $h_0+h_1$, where
\[
 h_0=\binom94+4\binom72+9=219,\qquad h_1=2\binom83+6\cdot6=148.
\]
For $h_0$ the three terms count monomials with zero, one, or two positive-weight factors and the same number of negative-weight factors; for $h_1$ they count monomials with one positive-weight factor and no negative one, or two positive and one negative. The limiting derivation has rank $\binom{13}{4}-367=348$, and semicontinuity of rank gives the claim.
\end{proof}

\begin{proposition}\label{prop:M10}
For $M=10$, a quartic in a nonempty Zariski open subset of $S_4$ has the avoidance property.
\end{proposition}

\begin{proof}
Here $q=55$ and $3q-3=162$. The projective family of matrices that are a scalar plus a matrix of rank at most one has dimension at most $1+(2M-1)-1=19$, and Lemma~\ref{lem:rank-bound} gives the rank bound $\binom{12}{3}=220>162+19$ on it. On the complement, Lemma~\ref{lem:rank348} gives $348>162+99$, where $99=\dim\PP(\End V)$. Both families are constructible, and Lemma~\ref{lem:incidence} applies to each stratum of a finite locally closed stratification.
\end{proof}

\subsubsection{Six variables}
Throughout this subsection $M=6$, so $\dim S_4=126$, $q=21$ and $\dim\mathcal Q_3=60$. The uniform rank bound is not enough on all of $\PP(\End V)$. We stratify matrices by their invariant factors and treat rank-one endomorphisms separately.

For $A\in\End V$, let $\lambda=(\lambda_1\geq\lambda_2\geq\cdots)$, $\sum_i\lambda_i=6$, be the degrees of its monic nonconstant invariant factors, and let $N_\lambda$ be the nilpotent matrix with Jordan blocks of lengths $\lambda_i$. Put $\lambda'_j=\#\{i:\lambda_i\geq j\}$ and $c_\lambda=\sum_j(\lambda'_j)^2$.

\begin{lemma}\label{lem:invariant-factors}
For matrices with invariant-factor degrees $\lambda$:
\begin{enumerate}
\item $\rank D_A\geq\rank D_{N_\lambda}$ on $S_4$;
\item the centralizer of $A$ in $\End V$ has dimension $c_\lambda$;
\item the subfamily of $\PP(\End V)$ on which $D_A:S_4\to S_4$ is singular has dimension at most
\begin{equation}\label{eq:plambda}
 p_\lambda=36-c_\lambda+\lambda_1-2.
\end{equation}
\end{enumerate}
\end{lemma}

\begin{proof}
Put $A$ in rational canonical form. For a monic invariant factor $p(z)$ of degree $e$, the polynomial $p_t(z)=t^ep(z/t)$ has coefficients polynomial in $t$ and specializes to $z^e$ at $t=0$; for $t\neq0$ its companion matrix is conjugate to $t$ times that of $p$. Doing this for all invariant factors gives a family whose general member is conjugate to $tA$ and whose special member is $N_\lambda$, and semicontinuity gives~(i).

View $V$ as a $\C[z]$-module. For monic $p_i,p_j$,
\[
 \dim\Hom_{\C[z]}\bigl(\C[z]/(p_i),\C[z]/(p_j)\bigr)=\deg\gcd(p_i,p_j),
\]
and since the invariant factors form a divisibility chain, the centralizer has dimension $\sum_{i,j}\min(\lambda_i,\lambda_j)=\sum_j(\lambda'_j)^2=c_\lambda$, proving~(ii).

For~(iii), list the invariant factors as $p_1\mid\dots\mid p_s$ and put $p_0=1$. The monic quotients $p_i/p_{i-1}$ have prescribed degrees summing to $\lambda_1$, and their coefficients parametrize all such chains by an affine space of dimension $\lambda_1$. On the direct sum of the companion matrices, $\det D_A$ is a polynomial in these coefficients, not identically zero: for $p_i(z)=(z-1)^{\deg p_i}$, $A$ is the identity plus a nilpotent and $D_A$ is $4\,\mathrm{id}$ plus a nilpotent. So the singular locus in the coefficient space has dimension at most $\lambda_1-1$; conjugation adds at most $36-c_\lambda$ dimensions, and projectivizing the scaling-invariant family subtracts one.
\end{proof}

For a Jordan block of length $e$ in $N_\lambda$ choose the $\mathfrak{sl}_2$-representation with weights $e-1,e-3,\dots,1-e$, and let $w_1,\dots,w_6$ be the combined list of weights. Put
\begin{equation}\label{eq:hj}
 h_j(\lambda)=\#\Bigl\{(b_1,\dots,b_6)\in\Z_{\geq0}^6:\ \sum_ib_i=4,\ \sum_iw_ib_i=j\Bigr\},
\end{equation}
the coefficient of $t^4z^j$ in $\prod_{i=1}^6(1-tz^{w_i})^{-1}$. The raising-operator argument of Lemma~\ref{lem:rank348} gives
\begin{equation}\label{eq:rank-formula}
 \rank D_{N_\lambda}=126-h_0(\lambda)-h_1(\lambda).
\end{equation}
For the nonminimal nonscalar types, the values obtained from~\eqref{eq:hj} are listed below.
\begin{center}
\renewcommand{\arraystretch}{1.1}
\begin{tabular}{cccccc}
\toprule
$\lambda$ & $c_\lambda$ & $h_0+h_1$ & $\rank D_{N_\lambda}$ & $p_\lambda$ & $\rank D_{N_\lambda}-60-p_\lambda$\\
\midrule
$(6)$ & 6 & 12 & 114 & 34 & 20\\
$(5,1)$ & 8 & 18 & 108 & 31 & 17\\
$(4,2)$ & 10 & 22 & 104 & 28 & 16\\
$(4,1,1)$ & 12 & 26 & 100 & 26 & 14\\
$(3,3)$ & 12 & 26 & 100 & 25 & 15\\
$(3,2,1)$ & 14 & 32 & 94 & 23 & 11\\
$(3,1,1,1)$ & 18 & 46 & 80 & 19 & 1\\
$(2,2,2)$ & 18 & 36 & 90 & 18 & 12\\
$(2,2,1,1)$ & 20 & 46 & 80 & 16 & 4\\
\bottomrule
\end{tabular}
\end{center}
For the smallest margin, $\lambda=(3,1,1,1)$, the weights are $2,0,-2,0,0,0$; no monomial has weight one, and $h_0=\binom73+\binom53+1=46$. Here $c_\lambda=18$ and $p_\lambda=19$, so the required inequality is $80>60+19$.

The table is a finite computation in the representation $\Sym^4(\C^6)$. The following script reproduces every row with integer arithmetic. Each monomial is represented by its nondecreasing list of four variable indices, its weight is the sum of the four corresponding weights, and the rank is given by~\eqref{eq:rank-formula}.
\begin{verbatim}
from itertools import combinations_with_replacement

partitions = [(6,), (5,1), (4,2), (4,1,1), (3,3),
              (3,2,1), (3,1,1,1), (2,2,2), (2,2,1,1)]

for lam in partitions:
    weights = [d-1-2*j for d in lam for j in range(d)]
    kernel = sum(
        sum(weights[i] for i in monomial) in (0, 1)
        for monomial in combinations_with_replacement(range(6), 4)
    )
    centralizer = sum(
        sum(d >= j for d in lam)**2
        for j in range(1, lam[0]+1)
    )
    dimension = 36 - centralizer + lam[0] - 2
    rank = 126 - kernel
    assert rank > 60 + dimension
    print(lam, centralizer, kernel, rank,
          dimension, rank - 60 - dimension)
\end{verbatim}

\begin{lemma}\label{lem:not-rank-one}
The quartics $F$ with $D_AF\in\mathcal C_3$ for some nonzero $A$ that is not a scalar plus a matrix of rank at most one have proper Zariski closure in $S_4$.
\end{lemma}

\begin{proof}
Where $D_A$ is invertible, apply Lemma~\ref{lem:incidence} with $r=126$ and $p\leq35$; the inequality is $126>60+35$. On the singular locus, the invariant-factor loci are constructible, and Lemma~\ref{lem:invariant-factors} with the positive margins in the table gives the conclusion on each stratum. The two partitions of six missing from the table are $(1^6)$ and $(2,1^4)$. The first consists of scalars. For the second, divisibility forces the four linear invariant factors to be $z-a$ and the quadratic one to be divisible by $z-a$, so $A=a\,\mathrm{id}+B$ with $\rank B=1$; conversely every nonscalar matrix of this form has these degrees. Both types are excluded by hypothesis, and the stratification is finite.
\end{proof}

\begin{lemma}\label{lem:scalar-rank-one}
If $A=a\,\mathrm{id}+B$ with $a\neq0$ and $\rank B\leq1$, then $\rank D_A\geq91$ on $S_4$.
\end{lemma}

\begin{proof}
If $B$ is nilpotent or zero, $D_A=4a\,\mathrm{id}+D_B$ is invertible. Otherwise $B$ is diagonalizable with eigenvalues $b,0,0,0,0,0$, $b\neq0$. Decomposing $S_4$ by the degree $j$ in the distinguished eigenvector gives eigenvalues $4a+jb$, $0\leq j\leq4$, with multiplicities $\binom{8-j}{4-j}$. The eigenvalue for $j=0$ is nonzero, and at most one other eigenvalue vanishes, with multiplicity at most $\binom73=35$. So $\rank D_A\geq126-35=91$.
\end{proof}

The projective family in Lemma~\ref{lem:scalar-rank-one} has dimension at most $11$, and $91>60+11$, so Lemma~\ref{lem:incidence} treats it, including the nonzero scalars. The remaining matrices have rank one, and their images consist of quartics divisible by a linear form, where $\mathcal C_3$ is smaller.

\begin{lemma}\label{lem:divisible}
Let $\ell\in V$ be nonzero, $\dim V=6$. Then $\mathcal C_3\cap\ell S_3$ is the image of an algebraic parameter space of dimension at most $34$.
\end{lemma}

\begin{proof}
Write an element of $\mathcal C_3$ as $ab-c^2$ with $a,b,c\in S_2$, and let $\overline S=S/(\ell)$, a polynomial ring in five variables. Divisibility by $\ell$ means
\begin{equation}\label{eq:mod-ell}
 \bar a\,\bar b=\bar c^{\,2}\quad\text{in }\overline S.
\end{equation}
By unique factorization every solution has one of the forms
\begin{equation}\label{eq:solutions}
 (\bar a,\bar b,\bar c)=(u^2h,v^2h,uvh),\ h\in\overline S_2,\ u,v\in\C;\qquad(\bar a,\bar b,\bar c)=(u^2,v^2,uv),\ u,v\in\overline S_1.
\end{equation}
Indeed, if $\bar a,\bar b\neq0$, dividing by their greatest common divisor leaves coprime factors with square product, each a square up to a constant, of even common degree; the common divisor therefore has degree two or zero, giving the two forms. Constants and the sign of $\bar c$ are absorbed into $u,v$. If $\bar a$ or $\bar b$ is zero then $\bar c=0$, which is covered by the first form, as is the zero triple. The first family has dimension at most $15+2-1=16$, because of the redundancy $(h,u,v)\mapsto(t^{-2}h,tu,tv)$, and the second at most $10$. Each triple in the quotient has $3\dim(\ell S_1)=18$ parameters of lifts to $(S_2)^3$, so the triples $(a,b,c)$ satisfying~\eqref{eq:mod-ell} form a space of dimension at most $34$, whose image under $(a,b,c)\mapsto ab-c^2$ is $\mathcal C_3\cap\ell S_3$.
\end{proof}

\begin{lemma}\label{lem:rank-one}
The quartics $F$ with $D_AF\in\mathcal C_3$ for some $A$ of rank one have proper Zariski closure in $S_4$.
\end{lemma}

\begin{proof}
A rank-one derivation is $D_A=\ell\partial_v$ with $\ell\in V\setminus0$, $v\in V^*\setminus0$; its image on $S_4$ is $\ell S_3$, of dimension $56$, and the projective rank-one locus is $\PP(V)\times\PP(V^*)$, of dimension $10$. For fixed $A$, Lemma~\ref{lem:divisible} gives at most $34$ parameters for the values in $\mathcal C_3\cap\operatorname{im}D_A$, and the fibre in $F$ over each value is an affine space of dimension $126-56=70$. The incidence has dimension at most $10+34+70=114<126$. This count works algebraically as $\ell$ varies: cover $\PP(V)$ by its six standard charts and choose a linear splitting of $S\to S/(\ell)$ on each.
\end{proof}

\begin{proposition}\label{prop:M6}
For $M=6$, a quartic in a nonempty Zariski open subset of $S_4$ has the avoidance property.
\end{proposition}

\begin{proof}
Lemma~\ref{lem:not-rank-one}, Lemma~\ref{lem:scalar-rank-one} with the incidence estimate, and Lemma~\ref{lem:rank-one} cover all nonzero matrices. Outside the union of the exceptional closures,~\eqref{eq:DA} gives the avoidance property.
\end{proof}

\begin{proof}[Proof of Theorem~\ref{thm:generic-quartic}]
For $k=1,2$ and $k\geq3$ use Propositions~\ref{prop:M6}, \ref{prop:M10} and~\ref{prop:M13}. The exceptional closures are invariant under scaling of $F$, so the open set descends to $\PP(S_4)$, and it meets $B_M$, since both are nonempty open subsets of an irreducible space; let $U_M$ be the intersection. Proposition~\ref{prop:avoid-recon} gives the conclusion for every $[F]\in U_M$ and every $[G]\in B_M$.
\end{proof}

\begin{remark}\label{rem:fermat}
Not every smooth quartic has the avoidance property: for $F=x_1^4+\dots+x_M^4$ the element $x_1^4$ lies in $(J_F)_4\cap\mathcal C_3$. This is why the argument gives a generic statement, which the deformation theorem then extends to all smooth quartics.
\end{remark}

\subsection{Categorical Torelli for quartics}\label{sec:quartic-torelli}

Fix $k\geq1$ and $M=4k+2$, and let $X=V(f)\subset\PP^{M-1}$ be a smooth quartic. We first recover the even Jacobian algebra from the enhanced category, which gives a generic categorical Torelli theorem, and then apply the deformation theorem.

\subsubsection{An intrinsic algebra from the Serre functor}
For a smooth proper dg category $\AA$, morphisms between endofunctors are taken between their perfect bimodules. If $S_\AA$ is the Serre functor, the \emph{Hochschild--Serre algebra} has components
\[
 \HS^{p,q}(\AA)=\Hom(\Id_\AA,S_\AA^q[p]),
\]
with multiplication given by composition after applying the intervening functor. An enhanced equivalence identifies Serre functors, and hence these components and their products \cite[Theorem~3.4]{LRZ26}.

Let $\rho$ be the internal degree shift in the matrix-factorization description of $\Ku(X)$; it satisfies $\rho^4\simeq[2]$, and the Serre functor is $\rho^{-M}[M]$ \cite[Proposition~2.2]{LRZ26}. Hence
\begin{equation}\label{eq:quartic-serre}
 S_X\simeq\rho^{-2}[2k+2],\qquad\rho^2\simeq S_X^{-1}[2k+2].
\end{equation}
So the square of the degree shift is expressed through intrinsic data. Put $T_X=S_X^{-1}[2k+2]$ and form the graded algebra
\begin{equation}\label{eq:AX}
 \mathcal H_X=\bigoplus_{a\geq0}\Hom(\Id,T_X^a),\qquad u\cdot v=T_X^a(v)\circ u\quad(u:\Id\to T_X^a,\ v:\Id\to T_X^b),
\end{equation}
using coherent powers of $T_X$. The isomorphism~\eqref{eq:quartic-serre} induces compatible isomorphisms of all powers and identifies this multiplication with composition in the factorization model.

\begin{proposition}\label{prop:even-jac-intrinsic}
The unital subalgebra of $\mathcal H_X$ generated by $(\mathcal H_X)_1$ is isomorphic, as a graded $\C$-algebra, to $R^{(2)}_f$, with $(R_f)_{2a}$ in degree $a$. Every enhanced Fourier--Mukai equivalence between two such Kuznetsov components induces an isomorphism of their even Jacobian algebras.
\end{proposition}

\begin{proof}
By~\eqref{eq:quartic-serre}, the degree-$a$ piece of~\eqref{eq:AX} is $\Hom(\Delta,\Delta(2a))$, and the sector computation \cite[Proposition~3.10]{LRZ26} gives
\begin{equation}\label{eq:quartic-sectors}
 (\mathcal H_X)_a\simeq\begin{cases}(R_f)_{2a},&a\neq2k+1,\\(R_f)_{4k+2}\oplus\C^3,&a=2k+1.\end{cases}
\end{equation}
The first summand is the identity sector. Its multiplication is that of the Jacobian algebra \cite[Theorem~3.5 and Remark~3.7]{LRZ26}, that is, the product of factorization kernels computed in \cite[Section~5.2]{BFK14}, and the identity sector is closed under products. Since $k\geq1$, degree one is not the exceptional degree $2k+1$, so $(\mathcal H_X)_1=(R_f)_2$ lies in the identity sector. Every monomial of even degree is a product of quadratic monomials, so the subalgebra generated by $(\mathcal H_X)_1$ is $\bigoplus_a(R_f)_{2a}$ with its usual product; this is also the case $\gcd(4,M)=2$ of \cite[Lemma~3.8 and Proposition~3.11]{LRZ26}. An enhanced equivalence conjugates the Serre functors, hence identifies the algebras $\mathcal H_X$, their degree-one parts, and the subalgebras they generate.
\end{proof}

\subsubsection{Generic and global categorical Torelli}
Let $B_{4,M}\subset\PP(S_4)$ be the locus of smooth quartics.

\begin{theorem}\label{thm:generic-cat-quartic}
There is a nonempty Zariski open subset $U_M\subset B_{4,M}$ such that, if $[f]\in U_M$ and $[g]\in B_{4,M}$, an enhanced Fourier--Mukai equivalence $\Ku(V(f))\simeq\Ku(V(g))$ implies that $V(f)$ and $V(g)$ are projectively isomorphic.
\end{theorem}

\begin{proof}
Take the open set of Theorem~\ref{thm:generic-quartic}. Proposition~\ref{prop:even-jac-intrinsic} turns the equivalence into a graded isomorphism $R^{(2)}_f\simeq R^{(2)}_g$, and Theorem~\ref{thm:generic-quartic} gives the projective isomorphism; the second quartic may be any smooth member.
\end{proof}

\begin{proof}[Proof of Theorem~\ref{thm:quartic}]
Write $X=V(f)$, $Y=V(g)$, and let $\Phi:\Ku(X)\simeq\Ku(Y)$ be the given equivalence. In the universal family over $B=B_{4,M}$, Theorem~\ref{thm:deg4} gives $H^1(Z,T_Z)\simeq\HH^2(\Ku(Z))$ and $\HH^{<0}(\Ku(Z))=0$ for every fibre $Z$; the equation family has surjective Kodaira--Spencer maps and perfect relative projection kernels compatible with base change (Lemma~\ref{lem:equation-family}). Proposition~\ref{prop:correspondence} gives a smooth irreducible variety $T$, a point $t$, smooth maps $p_1,p_2:T\to B$ through $[f]$ and $[g]$, and a family of equivalences. Both maps are dominant, so $T^\circ=p_1^{-1}(U_M)\cap p_2^{-1}(U_M)$ is dense and open, and by Theorem~\ref{thm:generic-cat-quartic} the two fibres are projectively isomorphic at every complex point of $T^\circ$. Let $c:B\to\mathcal M_{4,M}$ be the map to the coarse moduli space of smooth polarized quartics, which is separated \cite[Theorem~1.7 and Corollary~1.8]{Ben13}. The equalizer of $cp_1$ and $cp_2$ is closed and contains $T^\circ$, hence $t$. So $(X,\OO_X(1))\simeq(Y,\OO_Y(1))$, and since the complete hyperplane linear systems give the embeddings, the isomorphism is projective linear.
\end{proof}

\begin{corollary}\label{cor:nondivisible}
Let $3\leq d<M$ with $d\nmid M$. For any smooth degree-$d$ hypersurfaces $X,Y\subset\PP^{M-1}$, an enhanced Fourier--Mukai equivalence $\Ku(X)\simeq\Ku(Y)$ implies that $X$ and $Y$ are projectively isomorphic.
\end{corollary}

\begin{proof}
By \cite[Theorem~3.12]{LRZ26}, only $d=4$, $M=4k+2$ remains, and the condition $d<M$ gives $k\geq1$. Theorem~\ref{thm:quartic} applies.
\end{proof}
\section{Uniform finiteness in the Calabi--Yau case}\label{sec:finiteness}

In this section we prove Theorem~\ref{thm:finite}. Subsection~\ref{sec:fermat} shows that, near the Fermat hypersurface, only finitely many hypersurfaces have a given Veronese algebra. Subsection~\ref{sec:bound} turns this finiteness into a bound depending only on $d$ and $m$, using the fact that a quasi-finite morphism of finite type has fibres of bounded size. Theorem~\ref{thm:deformation} then carries the bound to every smooth hypersurface.

\subsection{Generic finiteness from the Fermat algebra}\label{sec:fermat}

In this subsection we prove generic finiteness for the isomorphism classes of the algebras $B_X$ of Proposition~\ref{prop:veronese-intrinsic}. We work on a slice through the Fermat hypersurface, transverse to the action of $\GL_M$, and show that the map from the slice and the changes of basis to the space of top forms of $B_X$ has injective differential on a dense open subset. Throughout,
\begin{equation}\label{eq:range5}
 d\geq3,\quad m\geq2,\quad(d,m)\neq(3,2),\quad M=md,\quad c=m(d-2)\geq3,\quad h=\binom{M+d-1}{d}-M^2.
\end{equation}

\subsubsection{The Fermat slice}
Set $\rho=d-2$, $\beta=\rho(1,\dots,1)$, and
\[
 \mathsf A=\{u\in\Z^M_{\geq0}:|u|=d,\ u_i\leq \rho\}.
\]
The monomials of degree $d$ not in $\mathsf A$ are the $M^2$ monomials $x_i^{d-1}x_j$ (with $j=i$ allowed); two exponents cannot both be at least $d-1$ since $d\geq3$, so $\#\mathsf A=h$. Let $W$ have basis $e_u$, $u\in\mathsf A$, and consider the slice
\begin{equation}\label{eq:slice}
 f_s=\frac1d\sum_{i=1}^Mx_i^d+\sum_{u\in\mathsf A}s_ux^u,\qquad s\in\mathbf A^h.
\end{equation}
On an open neighbourhood $S^0$ of $0$, the hypersurfaces are smooth, the classes $x^u$ form a basis of $(R_{f_s})_d$, and $x^\beta$ is a basis of its socle. Normalize the trace by $\tr_s(x^\beta)=1$ and put
\begin{equation}\label{eq:topform}
 T_s(z)=\frac1{c!}\tr_s\Bigl(\Bigl(\sum_{u\in\mathsf A}z_ux^u\Bigr)^c\Bigr),
\end{equation}
a regular family in $E=\Sym^c(W^*)$. At the Fermat point,
\begin{equation}\label{eq:T0}
 T_0=\sum_{\substack{\sum_u\nu_u=c\\ \sum_u\nu_uu=\beta}}\frac{z^\nu}{\prod_u\nu_u!}.
\end{equation}
The following is the form of Macaulay duality used here; see \cite{IK99}.

\begin{lemma}\label{lem:apolar}
The graded algebra $B_s=R^{(d)}_{f_s}$ is the apolar algebra of $T_s$. Two such algebras are isomorphic if and only if their top forms differ by a linear change of variables and a nonzero scalar.
\end{lemma}

\begin{proof}
The algebra is generated by $W$ and has perfect pairings in complementary degrees, so an element $a\in\Sym^jW$ is a relation exactly when $\tr_s(ab)=0$ for all $b\in\Sym^{c-j}W$, that is, when $a(\partial)T_s=0$. An algebra isomorphism is determined by its degree-one map and transports the one-dimensional socle functional up to a scalar. Conversely, a linear change of the top form transports its apolar ideal, and over $\C$ a nonzero scalar on the form is absorbed by a scalar matrix.
\end{proof}

\subsubsection{A distribution lemma}

\begin{lemma}\label{lem:distribution}
Let $k\geq1$, $s$, $b_1,\dots,b_M$, $\alpha_1,\dots,\alpha_M$ be nonnegative integers with $\sum_j\alpha_j=ks$ and $\alpha_j\leq kb_j$. Then $\alpha$ is a sum of $k$ integral vectors of degree $s$, each bounded coordinatewise by $b$.
\end{lemma}

\begin{proof}
Build the columns of a $k$-row array one at a time. In column $j$, put $\lfloor\alpha_j/k\rfloor$ in every row, and add one more in the $\alpha_j\bmod k$ rows with the smallest current sums. The row sums always differ by at most one, and the entries in column $j$ are at most $\lceil\alpha_j/k\rceil\leq b_j$. At the end the total is $ks$, so all row sums equal $s$.
\end{proof}

\begin{lemma}\label{lem:linear-f}
Assume $c\geq4$. Let $w=u-v$ for $u,v\in\mathsf A$, and set $\alpha=\beta+w$, $t=\sum_j(w_j)_+\leq d$, $b_j=\min(\rho,\alpha_j)$, and $\mathsf A_\alpha=\{u:|u|=d,\ 0\leq u_j\leq b_j\}$. If a function $f:\mathsf A_\alpha\to\C$ satisfies
\[
 \sum_{a=1}^cf(u_a)=0\quad\text{whenever }u_a\in\mathsf A_\alpha\text{ and }\sum_au_a=\alpha,
\]
then $f(u)=\ell\cdot u$ for some $\ell$ with $\ell\cdot\alpha=0$. Coordinates with $b_j=0$ may be discarded.
\end{lemma}

\begin{proof}
\emph{Step 1: differences along one exchange.} Choose active indices $p,q$, and set $C=b-e_p-e_q$ and $\gamma=\alpha-e_p-e_q$. For degree-$(d-1)$ vectors $T,T'\leq C$, put $E_j=\min\{C_j,\gamma_j-\max(T_j,T'_j)\}\geq0$. Such vectors $T$ exist, because $\sum_jC_j=cd-t-2\geq(c-1)d-2\geq d-1$. Since $\sum_jb_j=cd-t$ and $\sum_j(\gamma_j-C_j)=t$, we have $\sum_jE_j\geq(cd-2)-2(d-1)-t\geq(c-3)d\geq d-1$. Choose $T''\leq E$ of degree $d-1$. The vector $\alpha-T-T''-e_p-e_q$ is nonnegative, of degree $(c-2)d$, with coordinates at most $2\rho\leq(c-2)\rho$, so Lemma~\ref{lem:distribution} splits it into $c-2$ rows bounded by $\rho$; these rows are also bounded by $\alpha$, hence lie in $\mathsf A_\alpha$. Comparing the two decompositions whose remaining two rows are $(T+e_p,T''+e_q)$ and $(T+e_q,T''+e_p)$ shows that $f(T+e_p)-f(T+e_q)$ does not depend on $T$; using $T'$ gives the same for every choice. Call this difference $\delta_{pq}$.

\smallskip\noindent\emph{Step 2: linearity.} For three active indices there is a degree-$(d-1)$ vector bounded by $b-e_p-e_q-e_k$, since the total capacity is at least $(c-1)d-3\geq d-1$; so $\delta_{pq}+\delta_{qk}=\delta_{pk}$, and $\delta_{pq}=\ell_p-\ell_q$. The graph of bounded degree-$d$ vectors under unit exchanges is connected: if $u\neq v$, choose $p,q$ with $u_p>v_p$ and $u_q<v_q$; then $u-e_p+e_q$ respects the bounds and is closer to $v$. Hence $f(u)-\ell\cdot u$ is constant, and adding this constant divided by $d$ to every $\ell_j$ absorbs it. Finally, Lemma~\ref{lem:distribution} gives a decomposition of $\alpha$, and the hypothesis then gives $\ell\cdot\alpha=0$.
\end{proof}

\subsubsection{The stabilizer of the Fermat top form}\label{subsec:stabilizer}
For $H=(h_{uv})\in\End(W)$ write $D_H=\sum h_{uv}z_u\partial_{z_v}$.

\begin{lemma}\label{lem:stabilizer}
In the range~\eqref{eq:range5}, the kernel of $H\mapsto D_HT_0$ consists of the operators
\begin{equation}\label{eq:stab}
 D_a=\sum_{u\in\mathsf A}(a\cdot u)\,z_u\partial_{z_u},\qquad\sum_ia_i=0,
\end{equation}
and has dimension $M-1$.
\end{lemma}

\begin{proof}
\emph{The case $c\geq4$.} Decompose $H$ by weights $w=u-v$. A weight-$w$ contribution has output multidegree $\alpha=\beta+w$, so different weights do not cancel. For one weight, define $f$ on $\mathsf A_\alpha$ by $f(u)=h_{u,u-w}$ if $u-w\in\mathsf A$ and $f(u)=0$ otherwise. The coefficient of $z^\nu$ in $D_HT_0$ is $(\prod\nu_u!)^{-1}\sum_u\nu_uf(u)$, so its vanishing is the hypothesis of Lemma~\ref{lem:linear-f}. If $w\neq0$, choose $i$ with $w_i>0$. Then $f$ vanishes on the face $u_i=0$, whose capacity is $\sum_{j\neq i}b_j\geq cd-t-\rho\geq(c-2)d+2>d$. For active $p,q\neq i$, the capacity of this face after subtracting $e_p+e_q$ is at least $d-1$, and comparing two vectors in the face that differ by $e_p-e_q$ gives $\ell_p=\ell_q$. So the coefficients on all active coordinates other than $i$ agree, and evaluating at a degree-$d$ vector in the face shows that this common value is zero. Inactive coordinates have $\alpha_j=0$ and do not enter. Since $\ell\cdot\alpha=0$ and $\alpha_i=\rho+w_i>0$, also $\ell_i=0$, so every nonzero weight vanishes. For $w=0$ the lemma gives $h_{uu}=a\cdot u$ with $a\cdot\beta=\rho\sum a_i=0$; conversely these operators annihilate $T_0$, and unit exchanges show that $a$ is determined by the diagonal entries.

\smallskip\noindent\emph{The case $c=3$.} Then $d=m=3$ and $M=9$. Index the basis by three-element subsets $I\subset\{1,\dots,9\}$; the top form is the sum of the products $z_Iz_Jz_K$ over partitions into triples. For an off-diagonal entry $(U,V)$ with $|U\setminus V|\geq2$, choose a partition of the complement of $V$ such that $U\setminus V$ meets two blocks; replacing $V$ by $U$ gives a monomial in which $U$ is the only block whose removal leaves disjoint blocks, so only this entry contributes to it, and the entry vanishes. For $U=C\cup\{b\}$, $V=C\cup\{a\}$ with $|C|=2$, write the entry as $c_C$. Replacing one block of a partition gives $c_C+c_{C'}=0$ for disjoint $C,C'$ avoiding $a,b$; among the seven available indices, any two such sets have a third set disjoint from both, so all $c_C$ are equal, and the disjoint-pair relation makes them zero. The diagonal entries $\lambda_I$ satisfy $\lambda_I+\lambda_J+\lambda_K=0$ on every partition. Comparing the partitions $(C\cup\{a\},C'\cup\{b\},K)$ and $(C\cup\{b\},C'\cup\{a\},K)$ for disjoint two-element sets $C,C'$ outside $\{a,b\}$ shows that $\lambda_{C\cup\{a\}}-\lambda_{C\cup\{b\}}$ is independent of $C$; if two such sets are not disjoint, a third set disjoint from both lets one make the comparison twice. Choosing a two-element set disjoint from three given indices gives the additive cocycle identity, so these differences have the form $a_p-a_q$. Connectedness of the exchange graph of triples makes $\lambda_I-\sum_{i\in I}a_i$ constant; adding one third of this constant to each $a_i$ and using the partition condition gives $\lambda_I=\sum_{i\in I}a_i$ with $\sum a_i=0$, that is,~\eqref{eq:stab}, of dimension eight.
\end{proof}

\subsubsection{Geometric first variations}
For $u\in\mathsf A$ and $i$ with $u_i>0$, define
\[
 \alpha_{u,i}=\beta+de_i-u,\qquad P_{u,i}=\sum_{\substack{|\nu|=c\\ \sum_v\nu_vv=\alpha_{u,i}}}\frac{z^\nu}{\prod_v\nu_v!}.
\]
The first variation of~\eqref{eq:topform} is
\begin{equation}\label{eq:first-variation}
 \psi_u=\frac{\partial T_s}{\partial s_u}\Big|_{s=0}=-\sum_{i:u_i>0}u_iP_{u,i}.
\end{equation}
Indeed, to first order the Jacobian relation replaces $x_i^{d-1}$ by $-s_uu_ix^{u-e_i}$, and a monomial of socle degree acquires a nonzero trace after one replacement exactly when its exponent is $\alpha_{u,i}$; the factorials in~\eqref{eq:topform} give~\eqref{eq:first-variation}, including monomials with repeated factors.

\begin{lemma}\label{lem:transverse}
For every pair $(u,i)$ as above there is a coefficient functional $\lambda_{u,i}\in E^*$ with
\[
 \lambda_{u,i}(P_{v,j})=\delta_{(u,i),(v,j)},\qquad\lambda_{u,i}(D_HT_0)=0\ \text{for all }H.
\]
In particular, $L_0=\langle\psi_u:u\in\mathsf A\rangle$ has dimension $h$ and meets the orbit tangent space only in zero.
\end{lemma}

\begin{proof}
In our range $c\geq d$: for cubics this is $m\geq3$, and for $d\geq4$ it follows from $m(d-2)\geq2(d-2)\geq d$. Put $b'_i=d-u_i-1$ and $b'_j=\rho$ for $j\neq i$. Then $1\leq b'_i\leq \rho$, $(\alpha_{u,i})_i=d+b'_i-1\leq db'_i\leq cb'_i$ and $(\alpha_{u,i})_j=\rho-u_j\leq cb'_j$. Lemma~\ref{lem:distribution} gives a factorization of $\alpha_{u,i}$ into $c$ degree-$d$ vectors, each with $i$-th entry at most $d-u_i-1$. The corresponding monomial of $P_{u,i}$ cannot occur in an orbit tangent vector, since producing it by one replacement of factors would need a factor with $i$-th entry at least $d-u_i$. The multidegrees $\alpha_{u,i}$ are distinct: the unique coordinate that is at least $d$ determines $i$, and then $u$. Normalizing the coefficient functional of the chosen monomial gives $\lambda_{u,i}$.
\end{proof}

\begin{proposition}\label{prop:finite-fibres}
There is a nonempty open subset $S^{\mathrm{good}}\subset S^0$ such that
\begin{equation}\label{eq:injective}
 \End(W)\oplus T_sS^0\longrightarrow E,\qquad(H,v)\longmapsto D_HT_s+(dT)_s(v),
\end{equation}
is injective for every $s\in S^{\mathrm{good}}$. Consequently, for any fixed graded algebra $B$, the set of pairs
\begin{equation}\label{eq:finite-pairs}
 \{(s,[g]):s\in S^{\mathrm{good}},\ g:B_s\xrightarrow{\ \sim\ }B\ \text{graded}\}
\end{equation}
is finite, where $[g]$ identifies isomorphisms differing by a common scalar on their degree-one maps.
\end{proposition}

\begin{proof}
\emph{Step 1: a direct sum at the Fermat point.} Let $K=\ker(H\mapsto D_HT_0)$ and choose $G=\sum_{u\in\mathsf A}g_ux^u$ with every $g_u\neq0$. For $D_a\in K$ one has $D_aP_{u,i}=(da_i-a\cdot u)P_{u,i}$. Suppose that $D_a\psi_G\in L_0+\operatorname{im}(H\mapsto D_HT_0)$, where $\psi_G=\sum_ug_u\psi_u$. Applying the functionals of Lemma~\ref{lem:transverse} shows that, for each $u$, the numbers $da_i-a\cdot u$ agree for all $i$ in the support of $u$, so the $a_i$ agree on every such support. Since $M\geq2d$, every pair of indices lies in the support of a squarefree monomial of degree $d$, so all $a_i$ agree, and $\sum a_i=0$ gives $a=0$. Hence
\begin{equation}\label{eq:direct-sum}
 \operatorname{im}(H\mapsto D_HT_0)+K\psi_G+L_0
\end{equation}
is a direct sum of dimension $(h^2-M+1)+(M-1)+h=h^2+h$.

\smallskip\noindent\emph{Step 2: injectivity along a line.} Choose a complement of $K$ in $\End(W)$. Along $s=tG$, use the orbit columns for that complement, the columns $t^{-1}D_HT_{tG}$ for $H\in K$, and all parameter-derivative columns. These are regular at $t=0$, and their limits are the independent columns in~\eqref{eq:direct-sum}. So a maximal minor is nonzero for general $t\neq0$, and removing the rescaling does not change the rank there. This proves~\eqref{eq:injective} on a nonempty open set.

\smallskip\noindent\emph{Step 3: finiteness.} The map $\PGL(W)\times S^{\mathrm{good}}\to\PP(E)$ has injective differential: a projective kernel vector lifts to an equation whose right side is a scalar multiple of $T_s$, and subtracting a scalar matrix reduces it to~\eqref{eq:injective}. So every fibre is zero-dimensional, and, being of finite type, finite. Lemma~\ref{lem:apolar} identifies the fibre with the stated isomorphism data. This also shows that the projective graded automorphism group of each $B_s$ on the good locus is finite.
\end{proof}

The image of this slice contains representatives of a nonempty open part of the moduli of hypersurfaces. Indeed, the orbit tangent space at the Fermat equation is spanned by the $M^2$ monomials $x_i^{d-1}x_j$, which complement the monomials indexed by $\mathsf A$, so the map from $\GL_M$ times the slice to the space of degree-$d$ equations is dominant, and so is its restriction to any nonempty open subset of the slice.
\subsection{A uniform bound}\label{sec:bound}

Proposition~\ref{prop:finite-fibres} shows that, for each graded algebra $B$, only finitely many hypersurfaces on the good slice have Veronese algebra isomorphic to $B$. Theorem~\ref{thm:uniform-special} needs more: one integer that bounds these numbers for all $B$ at once. This follows from a general fact about quasi-finite morphisms.

\begin{lemma}\label{lem:bounded-fibres}
Let $f:X\to Y$ be a quasi-finite separated morphism of finite type between noetherian schemes. Then there is an integer $N$ such that every fibre of $f$ has at most $N$ points.
\end{lemma}

\begin{proof}
By Zariski's Main Theorem \cite[Tag~05K0]{Stacks}, $f$ factors as an open immersion $X\subset T$ followed by a finite morphism $\pi:T\to Y$. For $y\in Y$, the fibre $T_y$ is the spectrum of the finite $\kappa(y)$-algebra $(\pi_*\OO_T)_y\otimes\kappa(y)$, so $\#X_y\leq\#T_y\leq\dim_{\kappa(y)}\bigl((\pi_*\OO_T)_y\otimes\kappa(y)\bigr)$. Since $\pi_*\OO_T$ is coherent, the sets $Z_n=\{y:\dim_{\kappa(y)}((\pi_*\OO_T)_y\otimes\kappa(y))\geq n\}$ are closed by Nakayama's lemma and decrease with $n$. As $Y$ is noetherian, the chain stabilizes, and since every fibre dimension is finite, the stable value is empty. So the fibre dimensions, and hence the fibre cardinalities of $f$, are bounded.
\end{proof}

\begin{proposition}\label{prop:uniform}
In the range~\eqref{eq:range5} there is an integer $\Bd_{d,m}$ such that, for every graded algebra $B$, at most $\Bd_{d,m}$ isomorphism classes of hypersurfaces on the good slice have Veronese algebra isomorphic to $B$.
\end{proposition}

\begin{proof}
Consider the morphism $\Phi:\PGL(W)\times S^{\mathrm{good}}\to\PP(E)$, $(g,s)\mapsto[g\cdot T_s]$, from the proof of Proposition~\ref{prop:finite-fibres}. Its differential is injective at every point, so $\Phi$ is unramified and in particular locally quasi-finite; it is of finite type between varieties, hence quasi-finite and separated. Let $\Bd_{d,m}$ be the integer of Lemma~\ref{lem:bounded-fibres} for $\Phi$. By Lemma~\ref{lem:apolar}, a graded isomorphism $B_s\simeq B$ is the same as a projective change of variables carrying $T_s$ to the top form of $B$, so every hypersurface on the good slice with Veronese algebra isomorphic to $B$ gives a point of the fibre of $\Phi$ over the class of that top form, and distinct isomorphism classes of hypersurfaces give distinct points, since they have distinct $s$. Hence there are at most $\Bd_{d,m}$ such classes.
\end{proof}

\begin{proof}[Proof of Theorem~\ref{thm:finite}]
By the dominance established at the end of Subsection~\ref{sec:fermat}, the image of $\GL_M$ times the good slice contains a nonempty open subset of smooth equations; saturating it under coordinate changes gives an invariant open subset $U$ whose members all have representatives on the good slice. By Proposition~\ref{prop:veronese-intrinsic}, any set of categorical partners in $U$ has mutually isomorphic algebras $B_X$, so by Proposition~\ref{prop:uniform} it has at most $\Bd_{d,m}$ elements. Proposition~\ref{prop:cy-hyp} gives surjective maps $\kappa$ and vanishing negative Hochschild cohomology throughout the range~\eqref{eq:range5}; the relative projectors extend over the equation family; and the coarse moduli space of smooth polarized hypersurfaces is separated \cite[Theorem~1.7 and Corollary~1.8]{Ben13}, with equality meaning projective equivalence. Theorem~\ref{thm:uniform-special} extends the bound from $U$ to every smooth hypersurface. The lower bound holds because $X$ is its own partner.
\end{proof}

\begin{remark}\label{rem:explicit}
The bound $\Bd_{d,m}$ can be made explicit. By the strong Lefschetz property of the Fermat Jacobian ring \cite{Sta80,Wat87}, a general linear form $L$ satisfies $L^d:(R_{f_s})_d\to(R_{f_s})_{2d}$ injective for $s$ near the Fermat point, and then a graded isomorphism $B_s\to B$ is determined by the images of the $M$ elements $L^{d-1}y_i$, where $y_1=L,y_2,\dots,y_M$ are linear coordinates: the images of all other monomials are recovered recursively from the multiplication in $B$. The isomorphisms are therefore the solutions of a system of equations of bounded bidegree in $\PP^h\times\PP^{Mh-1}$, and the multihomogeneous B\'ezout theorem \cite[Theorem~1.11]{DKS13}, divided by the order $d^{M-1}M!$ of the monomial symmetry group of the Fermat hypersurface, gives an explicit value of $\Bd_{d,m}$ in terms of $d$ and $m$.
\end{remark}

\section{Categorical Torelli in large degree}\label{sec:large}

We combine the deformation theorem with Voisin's generic Torelli theorem. The first step is to recover the primitive Hodge structure of a very general hypersurface from its Kuznetsov component.

\begin{proposition}\label{prop:hodge}
Fix $d\geq3$ and $m\geq2$, and put $M=md$ and $n=M-2$. There is an individual very-general locus $B_{\mathrm{vg}}$ in the space $B$ of smooth degree-$d$ hypersurfaces in $\PP^{M-1}$ such that, if $X\in B_{\mathrm{vg}}$, $Y\in B$ and $\Ku(X)\simeq\Ku(Y)$ by an enhanced equivalence, then $H^n_{\prim}(X,\Q)\simeq H^n_{\prim}(Y,\Q)$ as rational Hodge structures. If $n$ is odd, one may take $B_{\mathrm{vg}}=B$.
\end{proposition}

\begin{proof}
An enhanced equivalence preserves topological K-theory, its Hodge structure, and the Euler pairing \cite[Lemma~5.2 and Proposition~5.4]{Per22}.

\smallskip\noindent\emph{Odd $n$.} By the Lefschetz theorem and Poincar\'e duality, $H^j(X,\Q)=0$ for odd $j\neq n$, and exceptional objects contribute no odd topological K-theory. So \cite[Proposition~5.23]{Per22} gives $K^{\mathrm{top}}_{-n}(\Ku(X))\otimes\Q\simeq H^n(X,\Q)=H^n_{\prim}(X,\Q)$, and the same for $Y$.

\smallskip\noindent\emph{Even $n=2a$.} Set $K_X=K^{\mathrm{top}}_0(\Ku(X))\otimes\Q$ and $P_X=H^n_{\prim}(X,\Q)(a)$. Under the rational Chern character, the even cohomology of $X$ is the direct sum of $P_X$ and the Tate space $U_X=\bigoplus_{j=0}^n\Q H^j$, $H=c_1(\OO_X(1))$, with the summand of degree $2j$ twisted by $j$. The exceptional collection has length $I=M-d$, and its Chern characters lie in $U_X$. A primitive class is annihilated by $H$, whereas the Chern characters of the line bundles and the Todd class of $X$ are polynomials in $H$, so by Riemann--Roch $P_X$ is Euler-orthogonal to the exceptional collection on both sides; in particular $P_X\subset K_X$. Hence
\begin{equation}\label{eq:KX}
 K_X=P_X\oplus A_X,\qquad A_X=K_X\cap U_X,\qquad\dim_\Q A_X=(n+1)-I=d-1.
\end{equation}
For the dimension: the Euler matrix of the exceptional collection is triangular with ones on the diagonal, so its classes are independent and their orthogonality conditions on $U_X$ have rank $I$. The same Riemann--Roch calculation gives
\begin{equation}\label{eq:euler-prim}
 \chi(P_X,A_X)=\chi(A_X,P_X)=0,\qquad\chi(u,v)=(-1)^a\int_Xu\smile v\quad(u,v\in P_X).
\end{equation}
The Euler form on the full rational topological K-theory is nondegenerate, being the Poincar\'e pairing composed with duality and multiplication by the invertible Todd class; the semiorthogonal decomposition and the unipotent Euler matrix of the exceptional collection then show that its restriction to $K_X$ is nondegenerate. By~\eqref{eq:KX} and~\eqref{eq:euler-prim},
\begin{equation}\label{eq:PX-orth}
 P_X=A_X^{\perp_\chi}={}^{\perp_\chi}\!A_X\quad\text{inside }K_X.
\end{equation}

For a very general hypersurface, the primitive rational Hodge structure is simple unless it is entirely of type $(a,a)$; the monodromy argument is recalled in \cite[p.~451]{Voi13}. The latter does not happen here, because the residue description gives $H^{m(d-1)-1,m-1}_{\prim}(X)\simeq(R_f)_0\simeq\C$, with $(m(d-1)-1)-(m-1)=m(d-2)>0$. So $P_X$ has no nonzero rational Hodge classes at a very general point. By \cite[Theorem~1.1 and Corollary~1.3]{CDK95}, the locus where a nonzero primitive rational class becomes a Hodge class is a countable union of proper closed algebraic subsets of $B$; let $B_{\mathrm{vg}}$ be its complement. For $X\in B_{\mathrm{vg}}$,~\eqref{eq:KX} gives $A_X=\Hdg(K_X)$, of dimension $d-1$.

Now let $\Ku(X)\simeq\Ku(Y)$ with $X\in B_{\mathrm{vg}}$. The induced Hodge isomorphism $K_X\simeq K_Y$ shows that $\dim\Hdg(K_Y)=d-1$; since $A_Y$ consists of Hodge classes and has this dimension, $A_Y=\Hdg(K_Y)$. So the equivalence maps $A_X$ to $A_Y$, hence their Euler orthogonals to each other by~\eqref{eq:PX-orth}. Removing the Tate twist proves the proposition.
\end{proof}

\begin{theorem}[{\cite[Theorems~0.7(1) and~0.11]{Voi22}}]\label{thm:voisin}
There is an integer $d_0\geq5$, independent of $N$, such that the following holds whenever $d\geq d_0$, $N\geq d$ and $d\mid N+1$: a very general smooth degree-$d$ hypersurface $X\subset\PP^N$ is determined among all smooth degree-$d$ hypersurfaces in $\PP^N$ by the rational Hodge structure $H^{N-1}_{\prim}(X,\Q)$.
\end{theorem}

Here generic Torelli is understood as in the introduction of \cite{Voi22}, with rational Hodge structures and no prescribed polarization. The uniformity in the ambient dimension is \cite[Theorem~0.11]{Voi22}; we have increased the threshold so that $d_0\geq5$. Numbered references are to arXiv:2004.09310v3.

\begin{proposition}\label{prop:conditional}
Let $d\geq3$ and $m\geq2$ with $(d,m)\neq(3,2)$, and put $M=md$. Suppose that a very general smooth hypersurface of degree $d$ in $\PP^{M-1}$ is determined, up to projective isomorphism, among all smooth hypersurfaces of degree $d$ in $\PP^{M-1}$, by the rational Hodge structure on its primitive middle cohomology. Then any two smooth hypersurfaces $X,Y\subset\PP^{M-1}$ of degree $d$ with equivalent enhanced Kuznetsov components are projectively isomorphic.
\end{proposition}

\begin{proof}
Intersect the very general locus of Proposition~\ref{prop:hodge} with the one in the hypothesis. If $X$ lies in the intersection and $Y$ is any smooth hypersurface of the same degree and dimension, an equivalence $\Ku(X)\simeq\Ku(Y)$ gives an isomorphism of primitive rational Hodge structures by Proposition~\ref{prop:hodge}, and the hypothesis shows that $X$ and $Y$ are projectively isomorphic. To pass to all hypersurfaces we use the deformation theorem. Since $(d,m)\neq(3,2)$, the maps $\kappa$ are surjective at every point of the equation family and negative Hochschild cohomology vanishes (Proposition~\ref{prop:cy-hyp}); the exceptional line bundles extend over the family with perfect projectors compatible with base change; and the coarse moduli space is separated \cite[Theorem~1.7 and Corollary~1.8]{Ben13}. Theorem~\ref{thm:uniform-special}, with bound one in its very-general form, extends the statement to every pair of smooth hypersurfaces.
\end{proof}

\begin{theorem}\label{thm:large-div}
There is an integer $d_0\geq5$, independent of $m$, such that for every $d\geq d_0$ and $m\geq2$, any two smooth degree-$d$ hypersurfaces $X,Y\subset\PP^{md-1}$ with equivalent enhanced Kuznetsov components are projectively isomorphic.
\end{theorem}

\begin{proof}
Let $d_0$ be as in Theorem~\ref{thm:voisin}, and let $d\geq d_0$ and $m\geq2$. Put $M=md$; then $M-1\geq d$ and $d$ divides $M$, so by Theorem~\ref{thm:voisin} a very general smooth hypersurface of degree $d$ in $\PP^{M-1}$ is determined, among all smooth hypersurfaces of degree $d$ in $\PP^{M-1}$, by its rational primitive Hodge structure. The resulting isomorphism is projective, since both Picard groups are generated by the hyperplane class, which an isomorphism preserves. Since $d_0\geq5$, we have $(d,m)\neq(3,2)$, and Proposition~\ref{prop:conditional} applies.
\end{proof}

\begin{remark}\label{rem:d0}
The threshold $d_0$ is uniform in $m$ but not explicit.
\end{remark}

\begin{proof}[Proof of Theorem~\ref{thm:large}]
Choose $d_0$ as in Theorem~\ref{thm:large-div}, and let $d\geq d_0$ and $M>d$. If $d\mid M$, write $M=md$ with $m\geq2$ and apply Theorem~\ref{thm:large-div}. If $d\nmid M$, apply \cite[Theorem~3.12]{LRZ26}, whose exceptional quartic series does not occur since $d_0\geq5$.
\end{proof}

\appendix
\section{Generic Torelli for projective hypersurfaces}\label{app:generic}

In this appendix we prove Theorem~\ref{thm:classical}. By the work of Donagi and of Cox and Green, only the quartics in $\PP^{4k+1}$ remain, and for these we combine the algebraic reconstruction of Subsections~\ref{sec:even-jac} and~\ref{sec:generic-quartic} with Donagi's symmetrizer construction and the standard passage from infinitesimal variations of Hodge structure to generic Torelli \cite{Don83,DG84,CDT87,CDK95,Voi22}. Apart from those two subsections, the appendix uses only classical Hodge theory; in particular it does not use Kuznetsov components.

There is one conceptual point. The Hodge structure of a single hypersurface does not see the multiplication of its Jacobian ring, whereas its infinitesimal variation does. We therefore first recover the even Jacobian algebra from the infinitesimal variation (Subsection~\ref{app:ivhs}), and then use algebraic Hodge loci to show that, for a very general quartic, an isomorphism of Hodge structures comes from an isomorphism of infinitesimal variations (Subsection~\ref{app:hodge-loci}). Subsection~\ref{app:donagi} combines the result with Donagi's theorem.

We keep the notation of Section~\ref{sec:quartics}: $\dim V=M$, $S=\Sym V$, $J_F$ is the Jacobian ideal of a quartic $F\in S_4$, $R_F=S/J_F$, and $R^{(2)}_F=\bigoplus_{a\geq0}(R_F)_{2a}$ is the even Jacobian algebra. Let $B_M\subset\PP(S_4)$ be the locus of smooth quartics and $U_M\subset B_M$ the nonempty open subset of Theorem~\ref{thm:generic-quartic}. A very general point is one outside a countable union of proper closed algebraic subsets.

\begin{theorem}\label{thm:quartic-generic-torelli}
Let $k\geq1$. If $X\subset\PP^{4k+1}$ is a very general smooth quartic, $Y\subset\PP^{4k+1}$ is any smooth quartic, and
\[
 H^{4k}_{\prim}(X,\Q)\simeq H^{4k}_{\prim}(Y,\Q)
\]
as rational Hodge structures, then $X$ and $Y$ are projectively isomorphic. In particular, the same holds for an isomorphism of polarized integral Hodge structures.
\end{theorem}

\subsection{From the infinitesimal variation to the even Jacobian algebra}\label{app:ivhs}
Put $M=4k+2$ with $k\geq1$, and $\sigma=2M=8k+4$, the socle degree of the Artinian Gorenstein ring $R_F$. For a smooth quartic $X_F\subset\PP^{M-1}$ let $\Lambda_F=H^{4k}_{\prim}(X_F,\Z)$ and $H_F=\Lambda_F\otimes\Q$, polarized by the primitive cup-product form with the usual sign convention. By Griffiths' residue description and the infinitesimal residue calculation \cite{Gri69,CG80}, see also \cite[Section~2]{Don83} and \cite[Theorem~0.5]{Voi22},
\begin{equation}\label{eq:griffiths}
 H^{3k-j,\,k+j}_{\prim}(X_F)\simeq(R_F)_{4j+2}\qquad(0\leq j\leq2k);
\end{equation}
after compatible normalizations of the residue maps, the infinitesimal variation of Hodge structure is the multiplication
\begin{equation}\label{eq:ivhs}
 (R_F)_4\otimes(R_F)_{4j+2}\longrightarrow(R_F)_{4j+6}\qquad(0\leq j<2k),
\end{equation}
where $(R_F)_4$ is the tangent space of the moduli space at $[X_F]$, and the polarization pairs the $j$-th and $(2k-j)$-th pieces through
\begin{equation}\label{eq:polarization}
 (R_F)_{4j+2}\otimes(R_F)_{8k+2-4j}\longrightarrow(R_F)_\sigma\simeq\C.
\end{equation}
The normalizations depend only on the indices, so they do not affect isomorphisms of these tensors. The Hodge pieces are $(R_F)_2,(R_F)_6,\dots,(R_F)_{8k+2}$; the space $(R_F)_4$ enters only as the tangent space.

\begin{lemma}[{\cite[Proposition~6.2 and Lemma~6.3]{Don83}}]\label{lem:symmetrizer}
For a smooth quartic $X_F\subset\PP^{4k+1}$ with $k\geq1$, the infinitesimal variation of its primitive Hodge structure determines $R^{(2)}_F$ up to graded algebra isomorphism, functorially under isomorphisms of infinitesimal variations.
\end{lemma}

\begin{proof}
This is the quartic case of Donagi's symmetrizer construction, in the form of \cite{DG84} for smooth hypersurfaces; see also \cite[Proposition~1.1 and the beginning of Section~2]{Voi22}. We check that the required inequalities hold. The lowest tensor in~\eqref{eq:ivhs} is
\[
 \beta:(R_F)_2\otimes(R_F)_4\longrightarrow(R_F)_6.
\]
The first symmetrizer recovers a second copy of $(R_F)_2$ inside $\Hom((R_F)_2,(R_F)_4)$, and the next one identifies the two copies up to scalar and hence recovers the kernel of $\Sym^2(R_F)_2\to(R_F)_4$. For these two steps the pairs of degrees are $(2,4)$ and $(2,2)$, and the inequalities of the cited lemma are
\[
 \max\{2,\sigma-4\}\geq3,\quad\sigma-6>0,\quad\max\{2,\sigma-2\}\geq3,\quad\sigma-4>0,
\]
which hold since $\sigma\geq12$. Since $J_F$ is generated by cubics, $(J_F)_{2a}=S_{2a-4}(J_F)_4$ for $a\geq2$, and the quadratic Veronese presentation of Lemma~\ref{lem:veronese} shows that the recovered quadratic kernel generates the ideal of relations of $R^{(2)}_F$. In \cite{Voi22} the polynomial structure is then recovered in large dimension; here it is recovered for all $k\geq1$ by Theorem~\ref{thm:generic-quartic}.
\end{proof}

\begin{corollary}\label{cor:ivhs-recon}
Let $[F]\in U_M$ and let $G$ define any smooth quartic in $\PP^{4k+1}$. If the primitive infinitesimal variations of Hodge structure of $X_F$ and $X_G$ are isomorphic, then $X_F$ and $X_G$ are projectively isomorphic.
\end{corollary}

\begin{proof}
An isomorphism of infinitesimal variations identifies the tangent spaces, the Hodge pieces, and the tensors~\eqref{eq:ivhs}, in particular the lowest tensors $(R_F)_2\otimes(R_F)_4\to(R_F)_6$ and $(R_G)_2\otimes(R_G)_4\to(R_G)_6$. By Lemma~\ref{lem:symmetrizer}, $R^{(2)}_F\simeq R^{(2)}_G$, and Theorem~\ref{thm:generic-quartic} gives the conclusion.
\end{proof}

\subsection{From Hodge structures to generic Torelli}\label{app:hodge-loci}
We now pass from isomorphisms of infinitesimal variations to isomorphisms of Hodge structures, following the strategy of Cox, Donagi and Tu \cite{CDT87}, with the algebraicity of Hodge loci from \cite{CDK95}; compare \cite[Section~1.1]{Voi22}.

\begin{proof}[Proof of Theorem~\ref{thm:quartic-generic-torelli}]
\emph{Step 1: the setting.} Let $\mathcal M=[B_M/\PGL(M)]$ be the moduli stack of smooth quartics, a smooth separated Deligne--Mumford stack \cite{Ben13}, and let $\mathcal H_\Z$ be its local system of primitive middle cohomology. We work in smooth \'etale charts of $\mathcal M$, or equivalently with the universal family over $B_M$; we replace $U_M$ by its saturation under $\PGL(M)$, which descends to an open substack. Write $m=\dim\mathcal M$.

\smallskip\noindent\emph{Step 2: the period map is an immersion.} The lowest piece of the differential of the period map is
\[
 (R_F)_4\longrightarrow\Hom\bigl((R_F)_2,(R_F)_6\bigr).
\]
If $\xi\in(R_F)_4$ annihilates $(R_F)_2$, it annihilates $(R_F)_{\sigma-4}=(R_F)_2(R_F)_{\sigma-6}$, and the perfect Gorenstein pairing between degrees $4$ and $\sigma-4$ gives $\xi=0$. So the period map is an immersion on every chart.

\smallskip\noindent\emph{Step 3: Hodge correspondences.} On $\mathcal M\times\mathcal M$ consider the polarizable integral variation of weight zero $\mathcal V=\mathcal Hom(\mathrm{pr}_1^*\mathcal H_\Z,\mathrm{pr}_2^*\mathcal H_\Z)$. A rational Hodge isomorphism $\phi$ gives an integral Hodge tensor $u=D\phi$ for some positive integer $D$. The theorem of Cattani, Deligne and Kaplan \cite[Theorem~1.1]{CDK95}, applied for every bound on the norm of $u$ with respect to a polarization of $\mathcal V$, gives countably many algebraic spaces of Hodge classes, finite over $\mathcal M\times\mathcal M$; they are constructed on $B_M\times B_M$ and descend, being equivariant. Invertibility of the flat tensor is locally constant on each of them, so we keep the components where it holds. Their irreducible components are countably many algebraic correspondences $Z_\alpha$, each labelled by its family of isomorphisms, and every pair of hypersurfaces with rationally isomorphic Hodge structures lies on one of them. Only the polarizability of $\mathcal V$ is used, so $\phi$ need not respect the polarizations.

\smallskip\noindent\emph{Step 4: dominant correspondences lie in the diagonal.} Let $p_1,p_2:Z_\alpha\to\mathcal M$ be the projections. On simply connected charts, the Hodge locus of a fixed $\phi$ is
\begin{equation}\label{eq:hodge-locus}
 \Phi_2(t_2)=\phi\,\Phi_1(t_1),
\end{equation}
where $\Phi_i$ are the marked period maps. By Step~2, if one tangent coordinate of~\eqref{eq:hodge-locus} vanishes, so does the other; hence the differentials of both projections are injective, both projections are locally quasi-finite, and $\dim Z_\alpha\leq m$. Suppose that $p_1$ is dominant. Then $\dim Z_\alpha=m$, so $p_2$ is dominant as well, and on a dense smooth open subset of $Z_\alpha$ both projections are \'etale and $p_1$ lands in $U_M$. Near every point of this subset the correspondence is a local biholomorphism between charts, and its labelled flat tensor is an isomorphism of the primitive variations; differentiating along the tangent isomorphism identifies the infinitesimal variations. By Corollary~\ref{cor:ivhs-recon}, the two quartics are projectively isomorphic at every such point. Let $\mathsf M$ be the separated coarse moduli space. The image of this dense subset lies in the diagonal $\Delta_{\mathsf M}\subset\mathsf M\times\mathsf M$, which is closed, so the image of all of $Z_\alpha$, including its singular points, lies in it.

\smallskip\noindent\emph{Step 5: the very general locus.} The closures of $p_1(Z_\alpha)$ for the correspondences on which $p_1$ is not dominant form a countable union of proper closed algebraic subsets. Let $X$ lie outside it. A rational Hodge isomorphism $H^{4k}_{\prim}(X,\Q)\simeq H^{4k}_{\prim}(Y,\Q)$ lies on a correspondence $Z_\alpha$ dominating the first factor, so by Step~4 the pair $([X],[Y])$ lies on $\Delta_{\mathsf M}$, and $X$ and $Y$ are projectively isomorphic. Pulling back to $B_M$ gives the very general locus of equations.
\end{proof}

\subsection{Combination with Donagi's theorem}\label{app:donagi}

\begin{proof}[Proof of Theorem~\ref{thm:classical}]
In the notation of the theorem, Donagi's generic Torelli theorem \cite[Section~6]{Don83} holds for $d\geq3$ and $M\geq4$, except for cubic surfaces and possibly for
\[
 d\mid M;\qquad d=4,\ M=4k+2\ (k\geq1);\qquad d=6,\ M=6k+3\ (k\geq1).
\]
Cox and Green removed the sextic exception \cite[Theorem, p.~121]{CG90}; their argument recovers the polynomial structure from the infinitesimal variation throughout the remaining range $d\nmid M$, except for the quartic series. Theorem~\ref{thm:quartic-generic-torelli} treats the quartic series for every $k\geq1$.

It remains to match the formulations. The version with rational Hodge structures and a very general first hypersurface is explained in \cite[Theorem~0.2, Remarks~0.1--0.4 and Section~1.1]{Voi22}, and follows from the argument of Subsection~\ref{app:hodge-loci}: infinitesimal Torelli for smooth hypersurfaces other than cubic surfaces makes a dominant Hodge correspondence generically an isomorphism of local variations, the reconstruction of the polynomial structure places it on the diagonal of the coarse moduli space, closure treats all its points, and excluding the countably many correspondences whose first projection is not dominant gives the very general hypersurface. In particular, the second hypersurface is arbitrary, and the same argument applies to the result of Cox and Green. For the polarized form of the reconstruction, a rational Hodge isomorphism suffices: for a very general hypersurface in this range the primitive polarization is unique up to a scalar \cite[Remark~0.1]{Voi22}, so the pullback of the second polarization is a multiple of the first, and after complexification the isomorphism can be rescaled. This proves Theorem~\ref{thm:classical}.
\end{proof}

\end{document}